\RequirePackage{fix-cm}
\documentclass[trsc,nonblindrev]{informs4}

\usepackage{lineno}
\modulolinenumbers[5]
\RequirePackage{tgtermes}
\RequirePackage{newtxtext}
\DeclareFontFamily{OMS}{ntxtlf}{}
\DeclareFontShape{OMS}{ntxtlf}{m}{n}{<->ssub*cmsy/m/n}{}
\RequirePackage{newtxmath}
\RequirePackage{bm}
\RequirePackage{endnotes}
\usepackage{natbib}
\bibpunct[, ]{(}{)}{,}{a}{}{,}
\def\bibfont{\small}
\def\bibsep{\smallskipamount}

\usepackage{float}
\usepackage{graphicx}
\usepackage{subfig}
\usepackage{amsmath}
\usepackage{amsfonts}
\usepackage[subnum]{cases}
\usepackage{booktabs}

\makeatletter
\newcommand*{\centerfloat}{%
  \parindent \z@
  \leftskip \z@ \@plus 1fil \@minus \marginparwidth
  \rightskip \leftskip
  \parfillskip \z@skip}
\makeatother

\OneAndAHalfSpacedXII
\EquationsNumberedThrough
\TheoremsNumberedThrough
\ECRepeatTheorems
\MANUSCRIPTNO{}

\usepackage{algorithm}
\usepackage[noend]{algpseudocode}
\algdef{SE}[DOWHILE]{Do}{doWhile}{\algorithmicdo}[1]{\algorithmicwhile\ #1}

\usepackage{array}
\usepackage{threeparttable}
\usepackage{multirow}
\usepackage{rotating}
\usepackage{makecell}
\usepackage{optidef}

\let\INFORMSproof\proof
\let\endINFORMSproof\endproof
\renewenvironment{proof}{\INFORMSproof{Proof.}}{\endINFORMSproof}
\usepackage{xcolor}
\usepackage{hyperref}
\hypersetup{
    colorlinks,
    linkcolor=red,
    citecolor=blue,
    urlcolor=blue
}

\begin{document}
\TITLE{Bulk Service Queueing for Transit Resilience under Short Random Service Suspensions}
\RUNAUTHOR{Mo et al.}
\RUNTITLE{Resilience of Public Transit Systems under Short Random Service Suspensions}

\ARTICLEAUTHORS{%
\begingroup
\setlength{\parindent}{0pt}
{\fontsize{11.5}{13.5}\selectfont\bfseries
\mbox{Baichuan Mo\textsuperscript{1}},
\mbox{Li Jin\textsuperscript{2,3,*}},
\mbox{Zhengzhong Ricky You\textsuperscript{1}},
\mbox{Zuo-Jun Max Shen\textsuperscript{4}},
\mbox{Haris N. Koutsopoulos\textsuperscript{5}},
\mbox{Jinhua Zhao\textsuperscript{6}}\par}
\vspace{3pt}
{\fontsize{8.5}{10.2}\selectfont
\textsuperscript{1}Department of Civil Engineering, Tsinghua University, Beijing 100084, China\\[-1pt]
\textsuperscript{2}UMich Joint Institute, Shanghai Jiao Tong University, Shanghai 200240, China\\[-1pt]
\textsuperscript{3}Tandon School of Engineering, New York University, Brooklyn, NY 11201, USA\\[-1pt]
\textsuperscript{4}Department of Industrial Engineering and Operations Research, University of California, Berkeley, Berkeley, CA 94720, USA\\[-1pt]
\textsuperscript{5}Department of Civil and Environmental Engineering, Northeastern University, Boston, MA 02115, USA\\[-1pt]
\textsuperscript{6}Department of Urban Studies and Planning, Massachusetts Institute of Technology, Cambridge, MA 02139, USA\\[-1pt]
\textsuperscript{*}Corresponding author: Li Jin, \href{mailto:li.jin@sjtu.edu.cn}{li.jin@sjtu.edu.cn}.}
\endgroup
}

\ABSTRACT{Short service suspensions are common in public transit systems, but their operational impacts remain difficult to quantify. We develop an analytical framework for measuring the resilience of a transit line under short random service suspensions. Vehicle movement is represented by a two state process in which vehicles either travel normally or stop during a suspension, and the induced stochastic headways enter a bulk service queueing model with finite vehicle capacity and passenger carryover. The model yields two classes of resilience indicators. Stability conditions determine whether station queues remain bounded, while closed form expressions characterize the mean and variance of station level queue length and waiting time. We construct an independent renewal approximation for headways whose common marginal distribution is obtained by taking the positive part of a raw headway formed from the incident adjusted scheduled headway and the difference between two independent compound Poisson exponential variables. The renewal approximation preserves the marginal effects of short suspensions while omitting serial dependence and delay propagation across multiple vehicles. Combining the resulting passenger arrival distribution with a Markov representation of passenger loads across stations allows the resilience indicators to be computed sequentially along the route. The stability condition identifies how incident frequency, incident duration, passenger demand, scheduled headway, vehicle capacity, and downstream available capacity jointly determine whether a station becomes unstable. Numerical experiments show that short suspensions disproportionately affect congested stations and that changes in incident duration and scheduled headway can dominate comparable changes in vehicle capacity. A recursive first in, first out simulation assesses the analytical approximations and clarifies the role of headway variability. The framework provides a tractable way to assess capacity, scheduled headway, and incident mitigation decisions for improving the resilience of transit lines exposed to frequent short disruptions.}

\maketitle

\section{Introduction}\label{sec_intro}
Public transit systems (PTSs) connect people with jobs, homes, and other daily activities and therefore play an important role in cities worldwide. Their service, however, is susceptible to frequent unplanned delays and disruptions. For example, \citet{mo2022impact} report an average of 75 incidents per day in the Chicago urban rail system, of which more than 75\% last less than 5 minutes. These short suspensions may result from signal failures, passenger behavior, or infrastructure problems. Quantifying these effects is therefore important for assessing transit resilience.

In this study, we use two operational resilience metrics for PTSs, \textbf{queue length} and \textbf{waiting time}. Specifically, we model a PTS as a bulk service queue and aim to derive closed form expressions for the mean and variance of passenger queue length and waiting time at a station under random service suspensions. To this end, we derive a stability criterion for passenger queues under a renewal representation of suspension induced headways. The representation captures throughput loss when a nonpositive raw headway is converted to bunching and the lost separation is not recovered in a later headway. We also characterize the steady state distribution of passenger queues, which naturally leads to the quantification of these metrics.

Queueing behavior at a public transit station is commonly represented by a bulk service queue, in which a vehicle serves waiting passengers up to its available capacity \citep{powell1981stochastic}. Most transit applications analyze one station, while \citet{islam2014bulk} extended the analysis to a route through a Markov representation of downstream passenger flows. These studies focus on normal operations. Models of bulk service under disruptions instead represent a server with operating and repair states and usually assume a fixed batch size \citep{madan1989single,krishnamoorthy2014queues}. That assumption does not represent a transit route where the available boarding capacity varies with the onboard load and passenger alighting. Mapping a server breakdown to route operations is also difficult because passenger flows link the service capacities of consecutive stations along the route.

To fill these research gaps, we develop a bulk service queueing framework that links suspension induced headways, stochastic passenger arrivals, station queues, and the vehicle capacity propagated from upstream stations. The route calculation allows available boarding capacity to vary with onboard load and passenger alighting. Tracking these interactions station by station shows how short disruptions erode transit resilience through queue growth, waiting time variability, and station instability. These measures can inform future control and planning strategies for public transit.

Within this framework, we explicitly model random service suspensions in a single route PTS. Such a system represents a bus route or one direction of a rail line and is a basic element of more complex public transit networks. Each vehicle may experience random suspensions while traveling along the route. Section \ref{sec_susp_speed} explains the operational meaning of this assumption, shows how it corresponds to many real world situations, and presents the assumption as a first step toward a general incident representation for PTSs.

Our contributions are summarized as follows.
\begin{enumerate}
\item \textit{We develop an analytical bulk service queueing framework for quantifying line level public transit resilience under short random service suspensions.} In particular, we represent service suspensions through a vehicle speed process with normal and disruption states. The speed process gives a direct operational meaning to the abstract notion of a server breakdown in public transit operations under short suspensions.

\item \textit{We develop an analytically tractable renewal approximation for headways under random service suspensions.} The common marginal distribution is constructed as the positive part of a raw headway formed from the incident adjusted scheduled headway and the difference between two independent compound Poisson exponential variables. Under the approximation, headways are independent across vehicle runs, and departure epochs accumulate from these nonnegative headways. As a result, the analytical service sequence prevents overtaking without imposing a pathwise first in, first out (FIFO) recursion. We then use a censored normal approximation with a point mass at zero to obtain a closed form moment generating function and derive the probability generating function (PGF) and moments of passenger arrivals within a headway. This is a new analytical contribution to bulk service queueing under short service suspensions.

\item \textit{We develop a route level probability propagation model for inter station passenger flows.} The model recursively links passenger alighting, available vehicle capacity, station queues, and departing vehicle loads. These links allow the resilience indicators to be computed sequentially at every station along the route.

\item \textit{We develop an interpolation based root solving algorithm that makes the analytical queueing indicators computationally implementable for realistic vehicle capacities.} The algorithm finds all complex roots required by the transform based bulk service queue model, addressing a numerical bottleneck that standard single start nonlinear solvers and existing root search heuristics may miss.

\end{enumerate}

The rest of this paper is organized as follows. Section \ref{sec_liter} reviews the literature on the bulk service queue problem, random service disruptions, and queueing models for PTSs. Section \ref{sec_model} presents the model settings for a single route system with random service suspensions. Section \ref{sec_analysis} presents the major analytical results and their derivations. Section \ref{sec_num_ex} provides numerical examples to illustrate the theoretical results and assesses the proposed approach using simulation. Section \ref{sec_conclusion} concludes the paper and discusses future research directions.

\section{Literature review}\label{sec_liter}
\subsection{Bulk service queues}
Bulk service queues serve customers in groups rather than individually. \citet{bailey1954queueing} formulated an early fixed capacity model with Poisson arrivals and an embedded Markov chain at service completion epochs. For transportation terminals, \citet{powell1981stochastic} developed transform methods for random passenger arrivals and bulk service. These models provide the queueing foundation for capacity constrained passenger boarding.

\subsection{Random service disruptions}\label{sec_liter_disrupt}
Models of queues with service breakdowns typically alternate between operating and repair states. \citet{madan1989single} introduced this structure to a fixed batch bulk service queue, while \citet{krishnamoorthy2014queues} review the broader breakdown literature. These models are not directly applicable to a transit route because the service capacity at a station depends on the onboard load propagated from upstream stations.

\subsection{Queueing models in public transit systems}
Transit queueing analyses are usually conducted at the station level. Under regular service, suppose that passengers arrive according to a Poisson process, can board the first arriving vehicle, and face a deterministic headway $H$. Their mean waiting time is then
\begin{align}
\mathbb{E}[W] = H/2,
\label{eq_mean_wt_1}
\end{align}
where $H$ is the service headway and $W$ is the passenger waiting time. With stochastic headways, \citet{osuna1972control} obtained
\begin{align}
\mathbb{E}[W] = \frac{1}{2} \cdot \left[\mathbb{E}[H] +\frac{\operatorname{Var}[H]}{\mathbb{E}[H]} \right],
\label{eq_mean_wt_2}
\end{align}
where $\mathbb{E}[H]$ and $\operatorname{Var}[H]$ are the expectation and variance of headways, respectively. When headways are deterministic, $\operatorname{Var}[H]=0$, and the model reduces to Eq. \ref{eq_mean_wt_1}.

Both Eqs. \ref{eq_mean_wt_1} and \ref{eq_mean_wt_2} assume that every passenger can board the first vehicle. Once capacity becomes binding, however, some passengers are left behind and must wait for a later vehicle \citep{mo2020capacity}. Transform based bulk service models account for this additional delay and provide closed form moments of queue length and waiting time under capacity constraints \citep{powell1985analysis}.

Extending these models from one station to a route requires the downstream propagation of passenger loads and available capacity. \citet{hickman2001analytic} characterized route level headway and passenger flow dynamics, and \citet{islam2014bulk} combined this structure with a bulk service model. Their formulation assumes Erlang headways and does not jointly propagate headway correlation and vehicle capacity. Empirical headway distributions can differ substantially from the Erlang distribution \citep{bellei2010transit}. Our model instead derives a suspension induced headway distribution and links it to downstream capacity propagation.

\subsection{Service interruptions in public transit systems}
Research on public transit interruptions includes impact analysis and operations control. Impact studies measure observed operational effects \citep{mo2022impact}. Control studies examine responses such as substitute bus services \citep{jin2016optimizing} and passenger path recommendations under demand or behavior uncertainty \citep{mo2023robust,mo2025individual}. The present study complements these approaches by deriving analytical resilience indicators, including stability conditions and station level moments of queue length and waiting time, within a route level capacity constrained model.

\section{Model}\label{sec_model}

This section defines the transit line, passenger flow primitives, and stochastic suspension process used in the analysis. We begin with the single route bulk service setting and describe how vehicles and passengers interact at each station. We then represent short service suspensions with a two state vehicle speed process and derive the stochastic headway used in the station level queueing analysis.

\subsection{Single route public transit system and vehicle movements}
\begin{figure}[htb]
\centering
\subfloat{\includegraphics[width=0.9\textwidth]{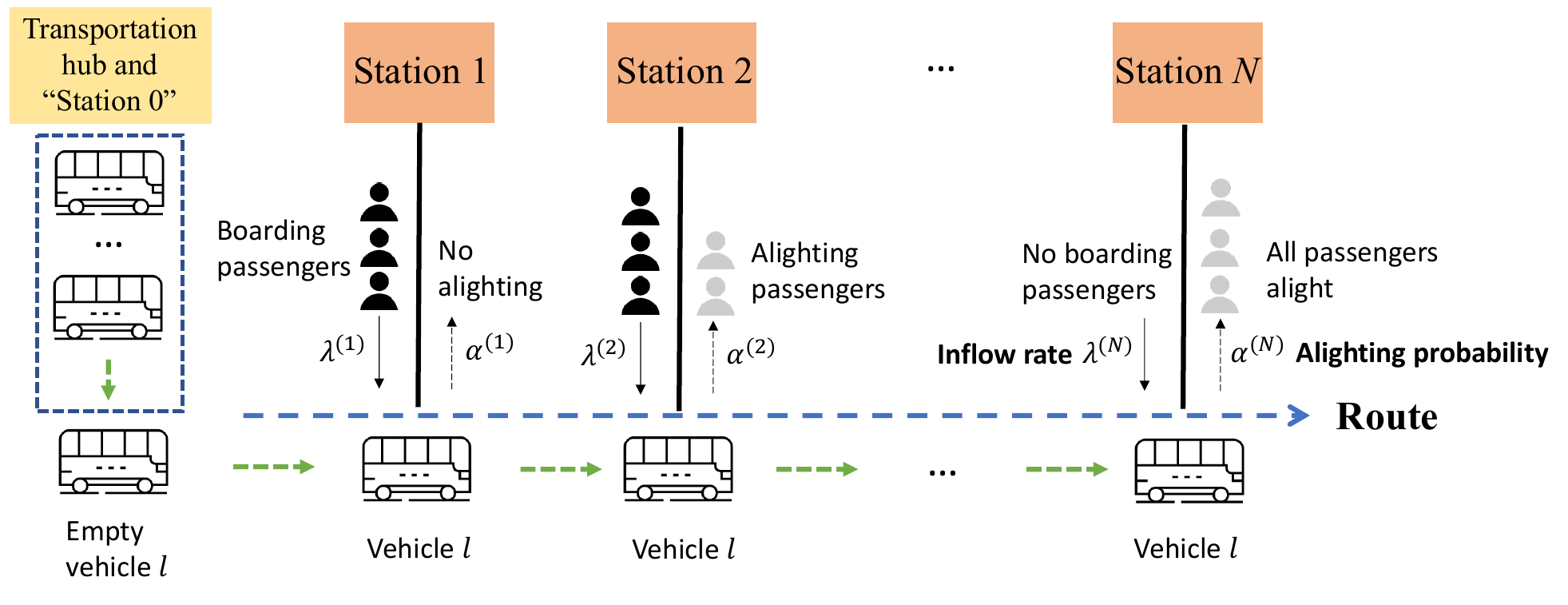}}
\caption{Schematic presentation of a single route public transit system}
\label{fig_net}
\end{figure}

Consider the single route PTS with $N$ stations shown in Figure \ref{fig_net}. Vehicles depart from a transportation hub, denoted by station 0, and visit stations 1 through $N$ in sequence. At station $n$, passenger arrivals follow a Poisson process with a fixed rate $\lambda^{(n)}$ during the period of interest. When a vehicle arrives, each onboard passenger independently alights with probability $\alpha^{(n)}$. Consequently, the number of alighting passengers at station $n$ follows a binomial distribution. Poisson arrivals and binomial alighting are common assumptions in the public transit literature \citep{hickman2001analytic}. Because the paper focuses on short service suspensions, we assume that passengers continue to wait and do not renege. The assumption is consistent with evidence that passenger abandonment becomes more important during substantially longer delays \citep{rahimi2019analysis}. Models of longer disruptions could relax the assumption by incorporating balking and reneging. Section \ref{append_model_assumptions} of the Electronic Companion (EC) appendix summarizes the modeling assumptions and their operational interpretations for the current setting.

\begin{figure}[htb]
\centering
\subfloat{\includegraphics[width=0.9\textwidth]{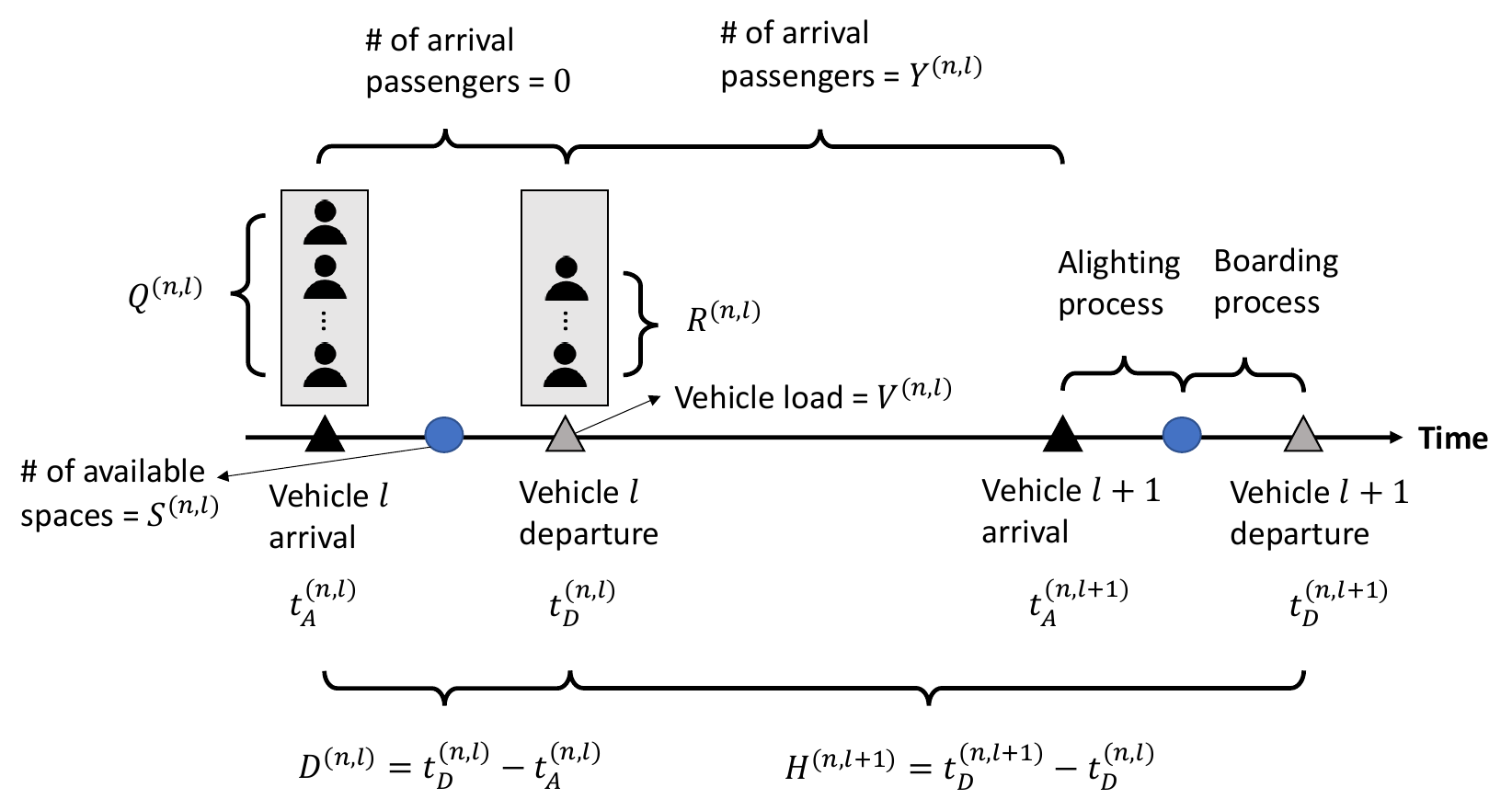}}
\caption{Vehicle and passenger interactions at station $n$ over time}
\label{fig_station}
\end{figure}

Let $l=1,2,\ldots$ denote the vehicle run index, with a smaller value indicating an earlier dispatch. Figure \ref{fig_station} summarizes the vehicle and passenger interactions at station $n$ over time. Let $t_A^{(n,l)}$ be the arrival time of vehicle $l$ at station $n$ and $t_D^{(n,l)}$ its departure time. For $l=2,3,\ldots$, we use $H^{(n,l)}$ provisionally for the service interval between vehicle runs $l-1$ and $l$. Section \ref{sec_headway} distinguishes its stationwise FIFO realization from the analytical renewal approximation. When a vehicle arrives, onboard passengers alight before queueing passengers board in their order of arrival. Let $Q^{(n,l)}$ be the number of queueing passengers when vehicle $l$ arrives at station $n$, $R^{(n,l)}$ the number left behind when it departs, and $Y^{(n,l)}$ the number arriving between $t_D^{(n,l)}$ and $t_A^{(n,l+1)}$. By definition,
\begin{align}
Q^{(n,l+1)} = R^{(n,l)} + Y^{(n,l)}.
\label{eq_relation_QRY}
\end{align}

Following \citet{powell1981stochastic}, we assume that dwell time is negligible relative to travel time and that no passenger arrives during dwell. Under this assumption, the number of passenger arrivals depends on the time between consecutive vehicles. For any service interval of length $h\geq0$, the conditional arrival count satisfies
\begin{align}
    Y^{(n,l)} \; |\; (\text{service interval}=h) \sim  \text{{\fontfamily{qcs}\selectfont Poi}}(\lambda^{(n)} h).
    \label{eq_Y_dist}
\end{align}
Thus, $Y^{(n,l)}$ counts passenger arrivals during the realized service interval.

To describe the same interaction from the vehicle side, let $S^{(n,l)}$ be the number of available spaces after passengers alight from vehicle $l$ at station $n$, and let $G^{(n,l)}$ be the number of passengers who remain onboard. By definition,
\begin{align}
G^{(n,l)} := C - S^{(n,l)},
\label{eq_Sn_Gn}
\end{align}
where $C$ is the vehicle capacity. Let $V^{(n,l)}$ denote the vehicle load when vehicle $l$ departs station $n$, which is also its load when it arrives at station $n+1$. Conditional on $V^{(n-1,l)}$, the number of passengers who alight from vehicle $l$ at station $n$ follows a binomial distribution.
\begin{align}
\left(V^{(n-1,l)} - G^{(n,l)}\right)\;|\;V^{(n-1,l)} \sim  \text{{\fontfamily{qcs}\selectfont Bin}}(V^{(n-1,l)}, \alpha^{(n)}).
\end{align}

\subsection{Random service suspensions and vehicle speed profile}\label{sec_susp_speed}
\begin{figure}[htb]
\captionsetup[subfigure]{justification=centering}
\centering
\subfloat[Vehicle speed profile]{\includegraphics[width=0.62\textwidth]{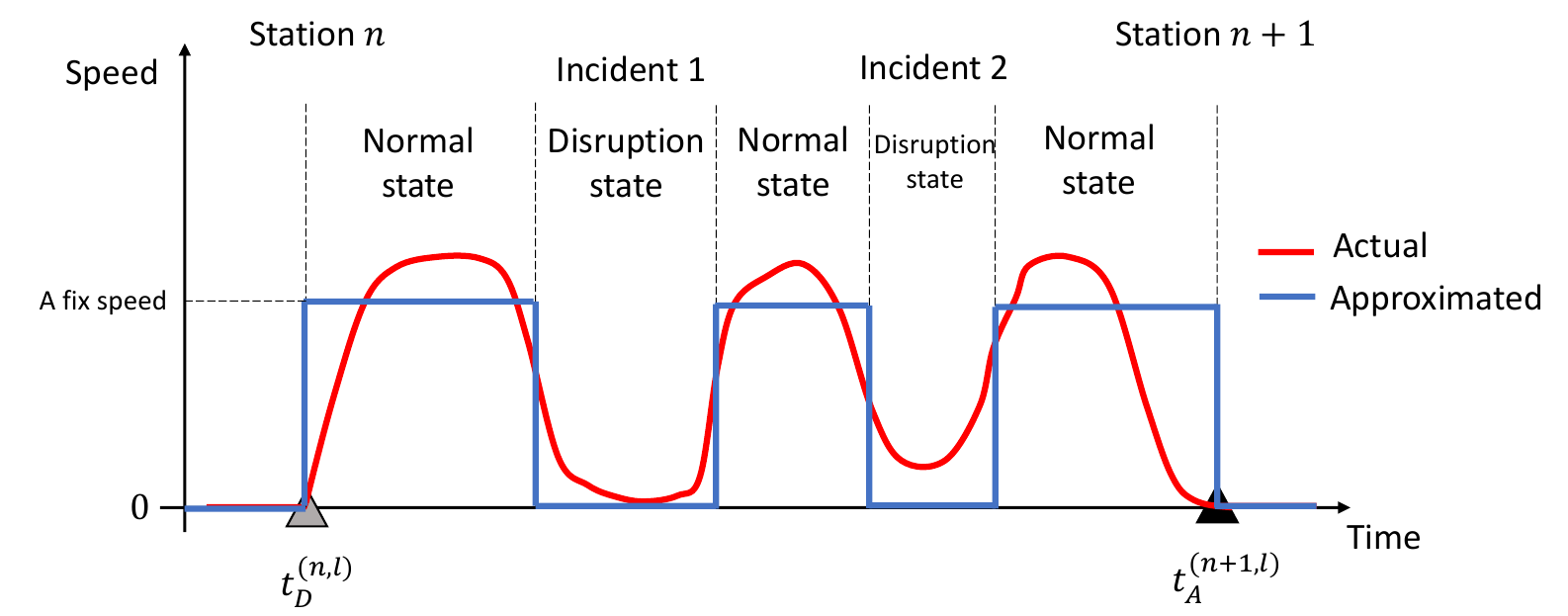}\label{fig_speed}}
\subfloat[Vehicle state transition]{\raisebox{0.045\textwidth}{\includegraphics[width=0.28\textwidth]{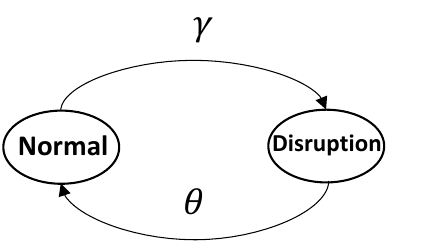}}\label{fig_failure_state}}
\caption{Representation of random service suspensions}
\end{figure}

We assume that random service suspensions may occur as a vehicle travels through the system. The red curve in Figure \ref{fig_speed} shows a possible speed profile for vehicle $l$ between stations $n$ and $n+1$ under such disturbances. Each incident slows or stops the vehicle. In bus systems, incidents may result from traffic congestion, accidents, vehicle engine problems, or driver and passenger behavior. In rail systems, they may result from signal or infrastructure failures and driver or passenger behavior. The speed profile provides a general representation of incidents, interruptions, suspensions, or disruptions that impede vehicle movement. We refer to the triggering event as an incident and to the resulting interruption of vehicle movement as a service suspension.

An actual vehicle speed profile under service suspensions can be complicated. For analytical tractability, we approximate it by the impulse function shown by the blue line in Figure \ref{fig_speed}. The impulse approximation separates the trajectory into traveling and stopping phases, which we call the normal and disruption states. In the normal state, a vehicle travels at a constant speed. When an incident occurs, the vehicle stops immediately and enters the disruption state. Over a sufficiently small interval $\Delta$, the probability of an incident is $\gamma\Delta$. Once the vehicle enters the disruption state, the incident duration follows an exponential distribution with rate $\theta$ and mean $1/\theta$. The vehicle state therefore follows the two state Markov process in Figure \ref{fig_failure_state} with state space $\{\text{Normal},\text{Disruption}\}$. Equivalently, the durations of the normal and disruption states are exponentially distributed with transition rates $\gamma$ for incident occurrence and $\theta$ for service recovery.

The impulse function provides a first step toward a more general incident representation. Different incident categories could be represented by normal and disruption states with heterogeneous occurrence rates and durations, producing a richer speed profile.

\subsection{Headway under random service suspensions}\label{sec_headway}
Under the assumed impulse function speed profile, all vehicles have the same fixed travel speed in the normal state. Consequently, when no incident occurs, all stations have the same deterministic headway $\overline{H}$, which we call the nominal scheduled headway. To relate this headway to route planning, let $\overline{E}$ denote the route cycle time required for a vehicle to travel from the transportation hub to the last station and return to the hub. Given route fleet size $\overline{F}$, the scheduled headway is $\overline{H} := \overline{E}/\overline{F}$.
Random service suspensions increase the expected route cycle time and therefore require a planning response. A transit agency can maintain the scheduled headway $\overline{H}$ by increasing the route fleet size, or it can maintain the same fleet size by increasing the scheduled headway. We consider the second response because it makes the effects of incidents on headway and service performance explicit, which is the focus of this paper.

Following this planning response, we assume that the transit agency has an estimate $\overline{D}$ of the average incident delay in the cycle time. Let $I^{(n,l)}$ be the total duration of all incidents while vehicle $l$ travels from the transportation hub to station $n$. Then $\mathbb{E}[I^{(N,l)}]$ is the expected incident duration for travel from the hub to the last station $N$. Assuming the same road conditions in both directions, the total estimated delay for the cycle trip is
\begin{align}
    \overline{D} &:=  2 \cdot \mathbb{E}[I^{(N,l)}]. \label{eq_e_delay}
\end{align}
Some transit agencies may plan for a delay above the mean, such as the 85th percentile. In that case, the cycle delay could be formulated as a more general function of the incident duration distribution. We use the expected delay in Eq. \ref{eq_e_delay} to keep the subsequent analysis transparent.
After incorporating this delay, the incident adjusted scheduled headway is $\overline{H}^{\text{Adj}} := (\overline{E}+\overline{D})/\overline{F} = \overline{H}+2\mathbb{E}[I^{(N,l)}]/\overline{F}$. The term $2\mathbb{E}[I^{(N,l)}]/\overline{F}$ is therefore the scheduled headway adjustment associated with incidents. Because $I^{(N,l)}$ is identically distributed across vehicle runs, the adjusted scheduled headway does not depend on the vehicle run index. The single route PTS then dispatches vehicles on time at intervals of $\overline{H}^{\text{Adj}}$.

Let $T^{(n)}$ be the travel time from the transportation hub to station $n$ when there is no incident. This travel time is fixed under the constant speed assumption. Since dwell time is negligible, the arrival and departure times at a station are identical in the model. Taking the dispatch time of vehicle 1 as time zero, the free departure time of vehicle $l$ at station $n$ is
\begin{align}
    \hat{t}_D^{(n,l)}
    := (l-1)\overline{H}^{\text{Adj}}+T^{(n)}+I^{(n,l)}. \label{eq_free_departure}
\end{align}
The free departure time in Eq. \ref{eq_free_departure} is the time that would result if vehicle order were unrestricted. Under this construction, a sufficiently long incident affecting vehicle $l-1$ can yield $\hat{t}_D^{(n,l)}<\hat{t}_D^{(n,l-1)}$. Vehicle $l$ would then appear to overtake vehicle $l-1$, although vehicle $l-1$ was dispatched first.

To rule out such reversals, we assume that vehicles serve each station in dispatch order and that overtaking is not allowed. The FIFO rule represents rail lines where passing is infeasible and bus services where passing is prohibited or uncommon at the route scale. It also keeps the vehicle index aligned with the service order observed by waiting passengers. If the free departure time of vehicle $l$ precedes the realized departure time of vehicle $l-1$, vehicle $l$ is held until vehicle $l-1$ departs. The realized departure times at station $n$ therefore satisfy
\begin{align}
    t_D^{(n,1)} := \hat{t}_D^{(n,1)}, \quad t_D^{(n,l)} := \max\left\{\hat{t}_D^{(n,l)},t_D^{(n,l-1)}\right\},
    \qquad \forall l\geq 2. \label{eq_fifo_departure}
\end{align}
Under Eq. \ref{eq_fifo_departure}, the stationwise FIFO headway is $H_{\mathrm{FIFO}}^{(n,l)}:=t_D^{(n,l)}-t_D^{(n,l-1)}=\max\{0,\hat{t}_D^{(n,l)}-t_D^{(n,l-1)}\}$. This headway depends on the realized departure time of the preceding vehicle and can therefore inherit delays from several earlier vehicles at the same station. Direct use of this recursion would therefore require the departure times of multiple vehicle runs to be analyzed jointly.

To avoid this recursive dependence, we impose the following approximation before constructing the analytical headway distribution.
\begin{assumption}[Independent renewal headways]\label{assumption_renewal}
At each station $n$, the analytical headways $\{H^{(n,l)}:l=2,3,\ldots\}$ are identically and independently distributed. Each headway is independent of all prior headways and of the queue and vehicle load states observed before the interval begins. The incident draws used to construct a headway are also independent of the exogenous passenger arrival process. The common marginal distribution is constructed from the free departure times of a generic pair of consecutive vehicles, with two new independent incident durations drawn for every headway in the renewal sequence.
\end{assumption}

Although the incident draws are independent of the exogenous Poisson passenger arrival process, the number of passengers arriving within an interval depends on its length through Eq. \ref{eq_Y_dist}. Separately, on time dispatch supports the assumed independence of the analytical headway from prior queue and vehicle load states. Independence across vehicle runs is an analytical approximation that is most defensible when the total incident delay of one vehicle rarely extends across several scheduled headways. The condition is consistent with service suspensions that are short relative to the scheduled headway. The approximation is also suitable for headway based operations in which holding or dispatch control resets the next service interval from the current departure \citep{daganzo2009headway}. It becomes less suitable when a long delay propagates through many subsequent vehicles. Section \ref{sec_num_ex} evaluates this limitation using a simulation that retains the recursive FIFO dynamics in Eq. \ref{eq_fifo_departure}.

Under Assumption \ref{assumption_renewal}, the common marginal distribution is obtained from the difference between the free departure times of a generic pair of consecutive vehicles. Accordingly, subtracting the two free departure times in Eq. \ref{eq_free_departure} gives the raw headway
\begin{align}
    \widehat{H}^{(n,l)} := \hat{t}_D^{(n,l)}-\hat{t}_D^{(n,l-1)} = \overline{H}^{\text{Adj}} + I^{(n,l)} - I^{(n,l-1)}, \qquad l=2,3,\ldots. \label{eq_headway_raw}
\end{align}
The deterministic term in Eq. \ref{eq_headway_raw} is $\overline{H}^{\text{Adj}}=\overline{H}+2\mathbb{E}[I^{(N,l)}]/\overline{F}$. For one marginal headway calculation, the two incident durations are independent. In a physical vehicle sequence, adjacent raw headways share one incident duration and are therefore dependent. The renewal model draws a new pair of incident durations for each analytical headway. If the difference between the two durations is sufficiently negative, $\widehat{H}^{(n,l)}$ becomes negative and cannot represent elapsed time between service events. We therefore define the analytical headway as
\begin{align}
    H^{(n,l)} := \max\{0,\widehat{H}^{(n,l)}\}, \qquad l=2,3,\ldots. \label{eq_renewal_headway}
\end{align}
Equation \ref{eq_renewal_headway} is not algebraically equivalent to the recursive FIFO rule in Eq. \ref{eq_fifo_departure}. A nonpositive raw headway becomes zero and represents bunching, but the next headway is redrawn instead of inheriting the accumulated delay. The approximation preserves vehicle order through nonnegative analytical service intervals, but it omits serial dependence and the propagation of accumulated delay to later vehicles.

\section{Analysis}\label{sec_analysis}

\begin{figure}[htb]
\centering
\subfloat{\includegraphics[width=0.9\textwidth]{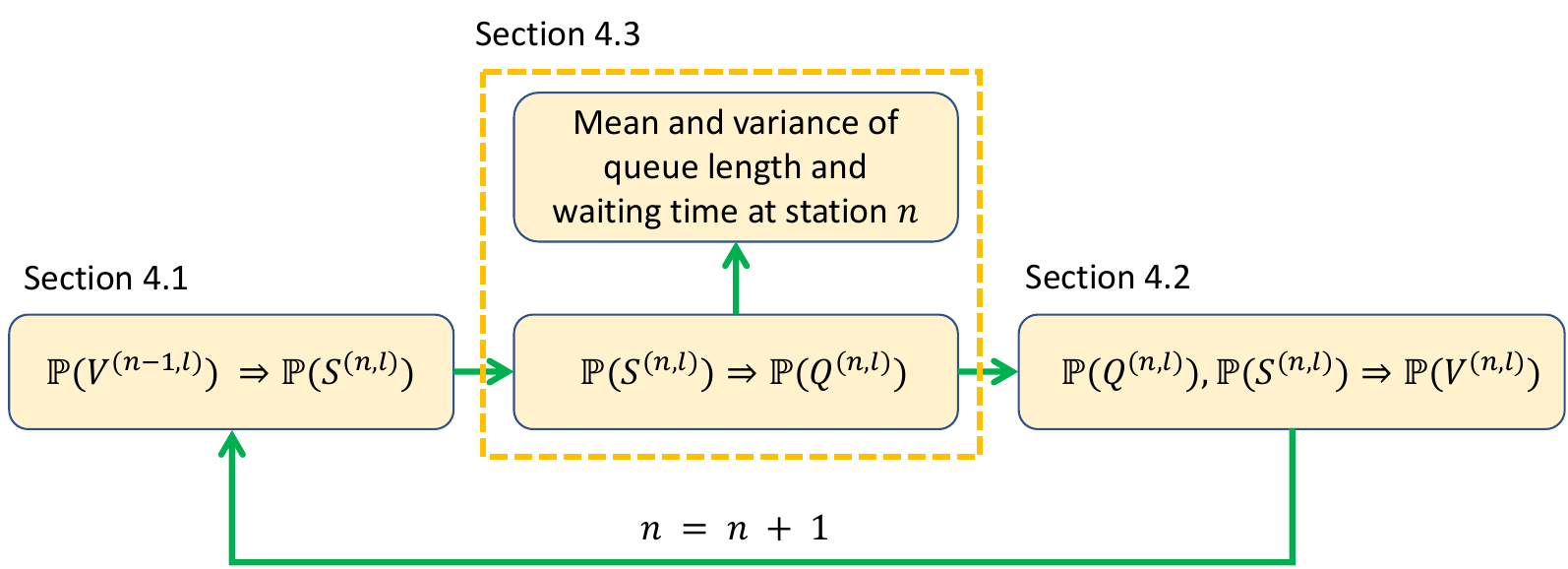}}
\caption{Analysis framework}
\label{fig_framework}
\end{figure}

The analysis derives the stability condition in Proposition \ref{prop_stability} and the mean and variance of passenger queue length and waiting time in Propositions \ref{prop_q_length} and \ref{prop_waiting}. Figure \ref{fig_framework} places these station level results within the route calculation by summarizing the propagation of the distributions of $S^{(n,l)}$, $V^{(n,l)}$, and $Q^{(n,l)}$. The calculation consists of three linked updates.
\begin{itemize}
\item Given the distribution of $V^{(n-1,l)}$, calculate the distribution of the available spaces $S^{(n,l)}$ after alighting at station $n$, as shown in Section \ref{sec_space_dist}.
\item Given the distribution of $S^{(n,l)}$, calculate the distribution of the queue $Q^{(n,l)}$ and the moments of queue length and waiting time, as shown in Section \ref{sec_que_anal}.
\item Given the distributions of $S^{(n,l)}$ and $Q^{(n,l)}$, calculate the distribution of the departing vehicle load $V^{(n,l)}$, as shown in Section \ref{sec_vehicle_load}.
\end{itemize}

Together, these components yield the distributions of $S^{(n,l)}$, $Q^{(n,l)}$, and $V^{(n,l)}$ for all $n=1,\ldots,N$. The recursion starts from $V^{(0,l)}=0$, which represents an empty vehicle departing the hub for the first station. The analysis focuses on steady state distributions as $l\to\infty$.

\subsection{Available vehicle space distribution}\label{sec_space_dist}
This subsection derives the steady state distribution of $S^{(n,l)}$ from that of $V^{(n-1,l)}$. Define $v_k^{(n,l)} := \mathbb{P}(V^{(n,l)}=k)$, $s_k^{(n,l)} := \mathbb{P}(S^{(n,l)}=k)$, and $g_k^{(n,l)} := \mathbb{P}(G^{(n,l)}=k)$ for $k=0,1,\ldots,C$. When the steady state probabilities exist, define $v_k^{(n)} := \lim_{l\to\infty}v_k^{(n,l)}$, $s_k^{(n)} := \lim_{l\to\infty}s_k^{(n,l)}$, and $g_k^{(n)} := \lim_{l\to\infty}g_k^{(n,l)}$. Let $V^{(n)}$, $S^{(n)}$, and $G^{(n)}$ denote random variables with these steady state distributions. Section \ref{sec_stable} establishes the stability condition under which the distributions exist.

\begin{proposition}[{Distribution of available vehicle space $S^{(n)}$}]\label{prop_space_dist}
For each $n=1,\ldots,N$, let $v^{(n-1)} := [v^{(n-1)}_0,\ldots,v^{(n-1)}_C] \in \mathbb{R}^{C+1}$ and $s^{(n)} := [s^{(n)}_0,\ldots,s^{(n)}_C] \in \mathbb{R}^{C+1}$. Let $g^{(n)} := [g^{(n)}_0,\ldots,g^{(n)}_C] \in \mathbb{R}^{C+1}$ denote the distribution of the vehicle load after alighting, with $g^{(n)}=v^{(n-1)}A^{(n)}$. The available space distribution is then given by the following reversal of the vehicle load index.
\begin{align}
s_k^{(n)} &= g_{C-k}^{(n)} \quad\quad \forall k = 0,1,\ldots,C. \label{eq_s}
\end{align}
The matrix $A^{(n)}$ has dimension $(C+1)\times(C+1)$. Its element $a_{ij}^{(n)}$ in row $i$ and column $j$ is the probability that $j$ passengers remain on board after alighting, conditional on $i$ passengers being on board before alighting. It is defined as
\begin{equation}
a_{ij}^{(n)} :=\left\{
\begin{aligned}
&\binom{i}{j}(1-\alpha^{(n)})^j(\alpha^{(n)})^{i-j}, & \quad\text{if $0\leq j\leq i\leq C$},\\
&0, & \quad\text{     otherwise}.
\end{aligned}
\right. \label{eq_A_mat}
\end{equation}
\end{proposition}

\subsection{Vehicle load distribution}\label{sec_vehicle_load}
This subsection states the steady state update for $V^{(n,l)}$ in terms of the distributions of $G^{(n,l)}$ and $Q^{(n,l)}$; Section \ref{sec_que_anal} subsequently derives the required queue distribution. Combining these two marginal distributions requires the following approximation for route propagation.

\begin{assumption}[Marginal route decoupling]\label{assumption_route_decoupling}
At each station $n$, the steady state remaining onboard load $G^{(n)}$ and the steady state local passenger queue $Q^{(n)}$ are independent. Therefore, $\mathbb{P}(G^{(n)}=i,Q^{(n)}=k)=\mathbb{P}(G^{(n)}=i)\mathbb{P}(Q^{(n)}=k)$ for all feasible $i$ and $k$.
\end{assumption}

Assumption \ref{assumption_route_decoupling} is exact at the first station because vehicles arrive from the hub without passengers and $G^{(1)}=0$. At a downstream station, the remaining onboard load results from boarding at upstream stations followed by binomial alighting, whereas the local queue results from the station specific Poisson arrival process and passengers previously left behind at that station. These quantities therefore arise from different passenger flow components. Independent renewal headways, alighting at intermediate stations, and repeated queue clearance under stable operation weaken the dependence induced by common headways and upstream capacity constraints. The approximation is most appropriate for short suspensions and stable stations, and it becomes weaker when long disruptions or persistent oversaturation jointly produce high onboard loads and long local queues. It permits the two marginal distributions to be propagated separately along the route. Section \ref{sec_num_ex} assesses the approximation through recursive simulation.

Define $q_k^{(n,l)} := \mathbb{P}(Q^{(n,l)}=k)$ and, when the steady state distribution exists, $q_k^{(n)} := \lim_{l\to\infty}q_k^{(n,l)}$. Let $Q^{(n)}$ denote a random variable with probabilities $\mathbb{P}(Q^{(n)}=k)=q_k^{(n)}$. The vector $q_{0:C-1}^{(n)} := [q_0^{(n)},\ldots,q_{C-1}^{(n)}]\in\mathbb{R}^{C}$ contains its first $C$ probabilities.

\begin{proposition}[{Distribution of vehicle load $V^{(n)}$}]\label{prop_load_dist}
Under Assumption \ref{assumption_route_decoupling}, the distribution of $V^{(n)}$ at each station $n=1,\ldots,N$, given $g^{(n)}$ and $q_{0:C-1}^{(n)}$, satisfies
\begin{align}
v^{(n)} = g^{(n)}B^{(n)},
\label{eq_v}
\end{align}
where $b_{ij}^{(n)}$, the element in row $i$ and column $j$ of $B^{(n)}$, is defined as follows.
\begin{equation}
b_{ij}^{(n)} :=\left\{
\begin{aligned}
&q^{(n)}_{j-i}, & \quad\text{     if $0\leq  i \leq  j < C$} \\
&1 - \sum_{k=0}^{C - i-1}  q^{(n)}_{k}, & \quad\text{     if $j = C$ and $0 \leq i < C$}\\
&1, & \quad\text{     if $i = j = C$}\\
&0, & \quad\text{     otherwise}
\end{aligned}
\right. \quad\quad i,j=0,1,\ldots,C. \label{eq_B_mat}
\end{equation}
\end{proposition}

\subsection{Queueing analysis at a station}\label{sec_que_anal}
Given the distribution $s^{(n)} := [s_0^{(n)},\ldots,s_C^{(n)}]\in\mathbb{R}^{C+1}$ of $S^{(n)}$, this subsection identifies the characteristic roots and determines $q_{0:C-1}^{(n)}$ and the mean and variance of passenger queue length and waiting time.
The transform calculations in this subsection are conditional on station stability, so the steady state queue and its PGF are well defined. Section \ref{sec_stable} establishes the stability condition independently from the transient queue recursion.
For compact notation, let $Y^{(n)}$ denote the steady state number of passenger arrivals within a headway at station $n$, and write $\overline{S}^{(n)}=\mathbb{E}[S^{(n)}]$. The quantities $\overline{\overline{S}}^{(n)}$ and $\overline{\overline{\overline{S}}}^{(n)}$ denote the second and third central moments of $S^{(n)}$, respectively, and the same convention applies to $Y^{(n)}$.

\begin{assumption}[Exogenous station service sequence]\label{assumption_station_service}
For each station $n$, once the marginal distribution $s^{(n)}$ has been propagated from upstream, the station level analytical model represents $\{S^{(n,l)}:l=1,2,\ldots\}$ as independent draws from $s^{(n)}$. This sequence is independent of the analytical headways, local passenger arrivals, and local queue history. The difference $Y^{(n,l)}-S^{(n,l+1)}$ has a finite mean and is nonconstant when $\overline{Y}^{(n)}=\overline{S}^{(n)}$.
\end{assumption}

\subsubsection{Probability generating function and characteristic roots.}\label{sec_pgf_Q}
We first derive the PGF of the steady state queue $Q^{(n)}$ and then identify the characteristic roots required by the queue probability calculation in Section \ref{sec_cal_q}. Throughout this station level derivation, $Q(z)$, $Y(z)$, $\textsc{Num}(z)$, and $\textsc{Den}(z)$ refer to station $n$, whose superscript is suppressed.

\begin{proposition}[{PGF of $Q^{(n)}$}]\label{prop_pgf_of_q}
Under Assumption \ref{assumption_station_service}, for each $n=1,\ldots,N$, the available space distribution $s^{(n)}$ yields the following PGF for the station queue $Q^{(n)}$.
\begin{align}
Q(z) = \frac{\sum_{u=0}^{C}s^{(n)}_u\left[\sum_{i=0}^u q_i^{(n)} (z^C-z^{C-u+i}) \right]}{\frac{z^C}{Y(z)} - \sum_{u=0}^{C}s_u^{(n)} z^{C-u}},
\label{eq_Qz1}
\end{align}
where $Y(z)$ is the PGF of $Y^{(n)}$.
\end{proposition}

Let $\textsc{Num}(z)$ and $\textsc{Den}(z)$ denote the numerator and denominator of $Q(z)$ in Eq. \ref{eq_Qz1}. Because $Q(z)$ is a PGF, it is analytic for $|z|<1$ and bounded for $|z|\leq1$. Every nonunit root of $\textsc{Den}(z)$ in the closed unit disk must therefore be canceled by a root of $\textsc{Num}(z)$ with the same multiplicity. The unit root is also a root of $\textsc{Num}(z)$ because direct substitution gives $\textsc{Num}(1)=0$. Since $\textsc{Num}(z)$ is a polynomial of degree at most $C$, the $C$ characteristic roots counted according to multiplicity determine its factorization up to a scale factor. Section \ref{sec_cal_q} uses this factorization to solve for $q_0^{(n)},\ldots,q_{C-1}^{(n)}$.

To locate the required roots, we multiply $\textsc{Den}(z)$ by $Y(z)$ to obtain the characteristic function
\begin{align}
\mathcal{F}^{(n)}(z) := z^C-Y(z)\sum_{u=0}^{C}s_u^{(n)}z^{C-u}.
\label{eq_characteristic_function}
\end{align}
The result of \citet{adan2006rouche} implies that $\mathcal{F}^{(n)}(z)$ has exactly $C$ roots in the closed unit disk when roots are counted according to multiplicity. Because $s_C^{(n)}>0$ and $Y(0)>0$, zero is not a root of $\mathcal{F}^{(n)}(z)$. At every nonzero root, $Y(z)$ is nonzero, so the same root and multiplicity apply to $\textsc{Den}(z)$. In particular, $z=1$ is a simple root because $[\mathcal{F}^{(n)}]'(1)=\overline{S}^{(n)}-\overline{Y}^{(n)}>0$. Interior roots need not be distinct \citep{oblakova2019exact}.
We therefore solve $\mathcal{F}^{(n)}(z)=0$ for the distinct root locations $\zeta_1,\ldots,\zeta_J$ in the closed unit disk, where $J$ denotes the number of distinct locations. Let $m_1,\ldots,m_J$ denote their corresponding multiplicities, which satisfy $\sum_{j=1}^{J}m_j=C$. For use in the queue calculation, we list the roots as $z_0^*,\ldots,z_{C-1}^*$ according to multiplicity and adopt the convention $z_0^*=1$. Section \ref{sec_cal_q} uses this list to determine $q_0^{(n)},\ldots,q_{C-1}^{(n)}$, while Section \ref{sec_root_solving} presents the numerical procedure for constructing it.

\subsubsection{Queue length distribution.}\label{sec_cal_q}
Section \ref{sec_pgf_Q} identifies the $C$ characteristic roots required to factor the numerator of $Q(z)$. We now translate that factorization into the queue probabilities $q_0^{(n)},\ldots,q_{C-1}^{(n)}$. First, the normalization condition $Q(1)=1$ determines the scale factor. Matching the resulting polynomial coefficients then yields the required probabilities.

\begin{proposition}[{Distribution of queue length $Q^{(n)}$}]\label{prop_solve_q}
For every station $n=1,\ldots,N$ with $s_C^{(n)}>0$, the distribution of $S^{(n)}$, the roots $z_0^*,\ldots,z_{C-1}^*$ listed according to multiplicity, and $\overline{Y}^{(n)}$ determine the probabilities $q_{0:C-1}^{(n)}$ through
\begin{equation}
q_0^{(n)} = \frac{1}{s_C^{(n)}} (\overline{S}^{(n)} - \overline{Y}^{(n)})\prod_{i=1}^{C-1}\frac{z_i^*}{z_i^*-1},
\qquad q_{0:C-1}^{(n)} = \widetilde{\eta}^{(n)} (\Lambda^{(n)})^{-1}.
\label{eq_q0_final}
\end{equation}
where the coefficients $\eta_j^{(n)}$ are defined by $\sum_{j=0}^{C}\eta_j^{(n)}z^j:=\prod_{i=0}^{C-1}(1-z/z_i^*)$. Let $\widetilde{\eta}^{(n)} := [s_C^{(n)}q_0^{(n)}\eta_0^{(n)},\allowbreak s_C^{(n)}q_0^{(n)}\eta_1^{(n)},\allowbreak\ldots,\allowbreak s_C^{(n)}q_0^{(n)}\eta_{C-1}^{(n)}]\in\mathbb{R}^{C}$ and
\begin{equation}
 \Lambda^{(n)} := \left[
    \begin{array}{ccccc}
     s_C^{(n)} & s_{C-1}^{(n)} & \cdots & s_2^{(n)} & s_1^{(n)} \\
     0 & s_C^{(n)} & \cdots & s_3^{(n)} & s_2^{(n)} \\
     \vdots & \ddots & \ddots & \vdots & \vdots \\
     0 & \cdots & 0 & s_C^{(n)} & s_{C-1}^{(n)} \\
     0 & \cdots & \cdots & 0 & s_C^{(n)}
    \end{array}\right] \in \mathbb{R}^{C\times C}.
\end{equation}
A root of multiplicity $m$ appears $m$ times in the defining product and in the root sums used below. Because the roots depend on station $n$, the coefficients $\eta_j^{(n)}$ carry the same station superscript.
\end{proposition}

The condition $s_C^{(n)}>0$ holds at the first station because every vehicle arrives empty. At a downstream station, it holds whenever there is a positive probability that all onboard passengers alight or that the vehicle arrives empty. In a degenerate case with $s_C^{(n)}=0$, the station level transform can instead use the largest $u$ for which $s_u^{(n)}>0$ as the effective maximum service capacity without changing the queue process.

\subsubsection{Analytical formulation of mean and variance of queue length and waiting time.}
Once $q_0^{(n)},\ldots,q_{C-1}^{(n)}$ have been obtained, the PGF $Q(z)$ is fully determined. We can therefore recover the expectation and variance of queue length at station $n$ from its first two derivatives as follows.
\begin{align}
\mathbb{E}[Q^{(n)}] &= \sum_{k = 0} ^ \infty k q_k^{(n)} =\left.\frac{dQ(z)}{d{z}}\right\vert_{z = 1},\\
\operatorname{Var}[Q^{(n)}] &= \mathbb{E}[(Q^{(n)})^2] - \mathbb{E}[Q^{(n)}]^2 =\left.\frac{d^2Q(z)}{d{z^2}}\right\vert_{z = 1} + \mathbb{E}[Q^{(n)}] - \mathbb{E}[Q^{(n)}]^2.
\end{align}

\begin{proposition}[{Mean and variance of queue length}]\label{prop_q_length}
For each $n=1,\ldots,N$, given the distribution of $S^{(n)}$ and the PGF $Y(z)$, the mean and variance of $Q^{(n)}$ are
\begin{align}
\mathbb{E}[Q^{(n)}] &= \frac{\overline{\overline{S}}^{(n)} + \overline{\overline{Y}}^{(n)}  + (\overline{S}^{(n)} - \overline{Y}^{(n)})[1+2(\overline{S}^{(n)} - C)] - (\overline{S}^{(n)} - \overline{Y}^{(n)})^2}{2(\overline{S}^{(n)} - \overline{Y}^{(n)})} + \sum_{i=1}^{C-1} \frac{1}{1-z_i^*}, \label{eq_EQ}\\
\operatorname{Var}[Q^{(n)}] &= \frac{1}{12(\overline{S}^{(n)} - \overline{Y}^{(n)})^2}\bigg[-4(\overline{\overline{\overline{S}}}^{(n)} - \overline{\overline{\overline{Y}}}^{(n)})(\overline{S}^{(n)} - \overline{Y}^{(n)}) + 3(\overline{\overline{S}}^{(n)} + \overline{\overline{Y}}^{(n)})^2 \nonumber \\
& - [6({\overline{\overline{S}}}^{(n)} - {\overline{\overline{Y}}}^{(n)}) - 1](\overline{S}^{(n)} - \overline{Y}^{(n)})^2    - (\overline{S}^{(n)} - \overline{Y}^{(n)})^4 \bigg] - \sum_{i=1}^{C-1} \frac{z_i^*}{(1-z_i^*)^2}.
\end{align}
\end{proposition}

We next express waiting time in terms of the queue length observed at an arbitrary time.

\begin{proposition}[{Mean and variance of waiting time}]\label{prop_waiting}
For each $n=1,\ldots,N$ with $\lambda^{(n)}>0$, given the distribution of $S^{(n)}$ and the PGF $Y(z)$, the mean and variance of waiting time $W^{(n)}$ are
\begin{align}
\mathbb{E}[W^{(n)}] &= \frac{\overline{Q}_t^{(n)} }{\lambda^{(n)}}, \label{eq_EW}\\
\operatorname{Var}[W^{(n)}] &=  \frac{\overline{\overline{Q}}_t^{(n)} - \overline{Q}_t^{(n)}}{(\lambda^{(n)})^2}, \label{eq_VarW}
\end{align}
where ${Q}_t^{(n)}$ is the queue length at an arbitrary time, whereas $Q^{(n)}$ is observed at a vehicle arrival. To make the distinction explicit, we write $\overline{Q}_t^{(n)}=\mathbb{E}[Q_t^{(n)}]$ and $\overline{\overline{Q}}_t^{(n)}=\operatorname{Var}[Q_t^{(n)}]$. The following expressions link the arbitrary time moments to their arrival epoch counterparts.
\begin{align}
\overline{Q}_t^{(n)} &= \mathbb{E}[Q^{(n)}] - \overline{Y}^{(n)} + \frac{1}{2}\left({\overline{\overline{Y}}^{(n)}}/{\overline{Y}^{(n)}} + \overline{Y}^{(n)} - 1\right),\\
\overline{\overline{Q}}_t^{(n)} &= \operatorname{Var}[Q^{(n)}] - \overline{\overline{Y}}^{(n)} + \frac{1}{12(\overline{Y}^{(n)})^2}\left[4\overline{Y}^{(n)} \overline{\overline{\overline{Y}}}^{(n)}+ 6(\overline{Y}^{(n)})^2\overline{\overline{Y}}^{(n)} - (\overline{Y}^{(n)})^2 + (\overline{Y}^{(n)})^4 - 3(\overline{\overline{Y}}^{(n)})^2 \right].
\end{align}
\end{proposition}

Under the Poisson arrival and first come, first served assumptions, Eq. \ref{eq_EW} applies Little's law, and Proposition \ref{prop_waiting} follows from \citet{powell1985analysis}.

\begin{remark}
The queue length and waiting time expressions are consistent with the general $M/G^{[S]}/1$ bulk service results of \citet{powell1985analysis}. The contribution here is the specification of $Y(z)$ induced by random service suspensions and its integration with route level capacity propagation.
\end{remark}

\subsubsection{Headway distribution.}\label{sec_hdw}
As defined in Eq. \ref{eq_renewal_headway}, the analytical headway is $H^{(n,l)} := \max\{0,\widehat{H}^{(n,l)}\}$. For the steady state analysis, we suppress the vehicle run index and write the common marginal variables as $H^{(n)}$ and $\widehat{H}^{(n)}$. We derive the distribution of $H^{(n)}$ in two steps. We first characterize the distribution of the raw headway $\widehat{H}^{(n)}$ and then take its positive part to obtain the distribution of $H^{(n)}$. This distribution is required in Propositions \ref{prop_solve_q} to \ref{prop_waiting} because Eq. \ref{eq_Y_dist} implies that $Y^{(n)}\mid H^{(n)}$ is Poisson with parameter $\lambda^{(n)}H^{(n)}$.

Equation \ref{eq_headway_raw} shows that the raw headway depends on the difference between two independent copies of the total incident duration for travel from the hub to station $n$. We therefore begin the derivation by characterizing the distribution of $I^{(n,l)}$.

\begin{proposition}[{Distribution of incident duration}]\label{prop_hdw_distribution}
The total incident duration $I^{(n,l)}$ for vehicle $l$ traveling from the transportation hub to station $n$ follows a compound Poisson exponential distribution with Poisson rate $\gamma T^{(n)}$ and exponential rate $\theta$. It has the representation
\begin{align}
    I^{(n,l)} := \sum_{i=1}^{K}X_i, \qquad X_i\sim\text{{\fontfamily{qcs}\selectfont Exp}}(\theta),\quad i=1,\ldots,K, \qquad K\sim\text{{\fontfamily{qcs}\selectfont Poi}}(\gamma T^{(n)}). \label{eq_I_sum}
\end{align}
\end{proposition}

The moment generating function (MGF) of a compound Poisson exponential variable is
\begin{align}
M_{I^{(n,l)}}(t) = \mathbb{E}[e^{tI^{(n,l)}}] =  e^{\gamma T^{(n)} (\frac{\theta}{\theta - t} - 1)},  \qquad t < \theta.
\end{align}
Similarly, the MGF of $-I^{(n,l-1)}$ is 
\begin{align}
M_{-I^{(n,l-1)}}(t) = \mathbb{E}[e^{-tI^{(n,l-1)}}] =  e^{\gamma T^{(n)} (\frac{\theta}{\theta + t} - 1)},  \qquad t > -\theta.
\end{align}
The MGF of $I^{(n,l)}$ gives $\mathbb{E}[I^{(N,l)}] = \gamma T^{(N)}/\theta$. Using this expectation in Eq. \ref{eq_headway_raw} and the independence of the two incident duration copies gives the raw headway MGF in Proposition \ref{prop_headway_distribution1}.

\begin{proposition}[{MGF of raw headway}]\label{prop_headway_distribution1}
Under the setting of this study, for each $n=1,\ldots,N$, the MGF of $\widehat{H}^{(n)}$ can be expressed, for $|t|<\theta$, as
\begin{align}
 M_{\widehat{H}^{(n)}}(t) =e^{t(\overline{H} + \frac{2 \gamma T^{(N)}}{\theta\overline{F}})}  e^{\gamma T^{(n)}(\frac{2t^2}{\theta^2 - t^2})}.
\label{eq_mgf_H}
\end{align}
\end{proposition}

From the MGF of $\widehat{H}^{(n)}$, we obtain the mean and variance of the raw headway as
\begin{align}
    \mathbb{E}[\widehat{H}^{(n)}] &=  \overline{H} + \frac{2 \gamma T^{(N)}}{\theta\overline{F}},\\
    \operatorname{Var}[\widehat{H}^{(n)}] &= \frac{4 T^{(n)} \gamma}{\theta^{2}}.
\end{align}

\begin{remark}
The results show that random suspensions can increase the mean and variance of the raw headway. The impact on its mean occurs through the cycle time adjustment at the planning stage. Its variance increases with a higher incident rate $\gamma$ and a higher average incident duration $\frac{1}{\theta}$. The station specific marginal distributions also capture the downstream growth in headway variability described by \citet{hickman2001analytic}. In particular, $\operatorname{Var}[\widehat{H}^{(n)}]$ increases with the station index $n$ because $T^{(n)}$ increases.
\end{remark}

Having characterized the raw headway, the second step applies the positive part transformation in Eq. \ref{eq_renewal_headway}. A direct MGF representation of the resulting headway involves nested integrals over the two compound Poisson distributions. Although this representation is exact, it does not reduce to a tractable finite expression for general parameter values. Moreover, the root calculations would require repeated numerical integration or infinite series evaluation at every complex argument. To avoid this computational burden, we first approximate the raw headway $\widehat{H}^{(n)}$ by a normal distribution with the same mean and variance and then take its positive part.

\begin{figure}[htb]
\centering
\subfloat[Example 1]{\includegraphics[width=0.33\linewidth]{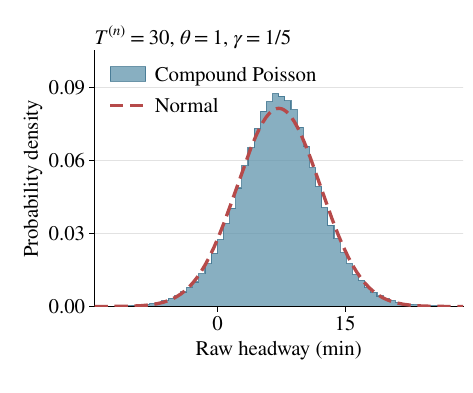}}
\subfloat[Example 2]{\includegraphics[width=0.33\linewidth]{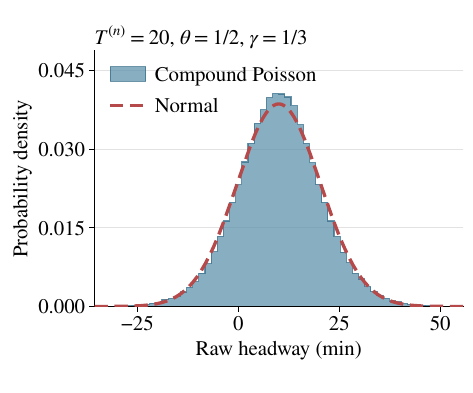}}
\subfloat[Example 3]{\includegraphics[width=0.33\linewidth]{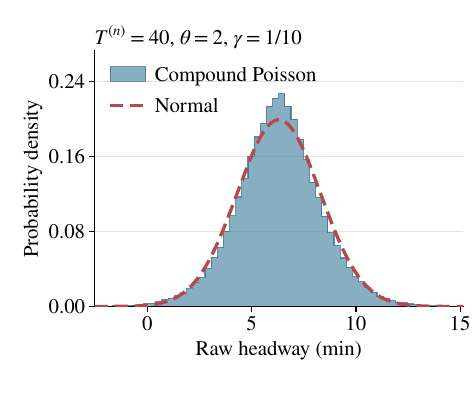}}
\caption{Numerical assessment of the normal approximation for the raw headway distribution}
\label{fig_emp_valid}
\end{figure}

We use the normal approximation for two reasons. First, $I^{(n,l)}$ is a random sum of independent exponential durations. A normal approximation becomes more accurate as the expected incident count increases, consistent with the setting of frequent short disturbances. A related central limit approximation for headway disturbances appears in \citet{daganzo2009headway}. Second, let $\widehat{H}^{(n)}_{\text{Normal}}$ denote a normal variable with the same mean and variance as $\widehat{H}^{(n)}$. Both variables are symmetric and have a third central moment of zero. Moreover, the MGF of $\widehat{H}^{(n)}_{\text{Normal}}$ is
\begin{align}
    M_{\widehat{H}^{(n)}_{\text{Normal}}}(t) =  e^{t(\overline{H} + \frac{2 \gamma T^{(N)}}{\theta\overline{F}})}  e^{\gamma T^{(n)}(\frac{2t^2}{\theta^2})},
\end{align}
which has the same second order expansion around $t=0$ as the exact MGF in Eq. \ref{eq_mgf_H}. The fourth central moment is $48(T^{(n)}\gamma)^2/\theta^4$ for $\widehat{H}^{(n)}_{\text{Normal}}$ and $[48(T^{(n)}\gamma)^2+48T^{(n)}\gamma]/\theta^4$ for $\widehat{H}^{(n)}$. The exact raw headway may therefore have heavier tails. Figure \ref{fig_emp_valid} numerically compares the two distributions for several values of $T^{(n)}$, $\theta$, and $\gamma$. The histograms are generated by sampling from the associated exponential and Poisson distributions. The normal distribution provides a close approximation in these examples while preserving the expected difference in tail behavior.

We next take the positive part of $\widehat{H}^{(n)}_{\text{Normal}}$. Denote the resulting headway approximation by $H^{(n)}_{\text{Normal}}$. It has support $[0,+\infty)$ and a probability mass at zero.

\begin{proposition}[{MGF of approximated headway}]\label{prop_headway_distribution}
For $\gamma>0$ and $n=1,\ldots,N$, the MGF of $H^{(n)}_{\text{Normal}}$ is given in Eq. \ref{eq_mgf_approximated_headway}.
\begin{equation}
\begin{aligned}
  M_{H^{(n)}_{\text{Normal}}}(t)
  &=\Phi\left(\frac{-(\overline{H}\theta + \frac{2\gamma T^{(N)}}{\overline{F}})}{2\sqrt{T^{(n)}\gamma}}\right) \\
  &\quad + e^{t(\overline{H} + \frac{2 \gamma T^{(N)}}{\theta\overline{F}})} e^{\gamma T^{(n)}(\frac{2t^2}{\theta^2})}
  \left[1-\Phi\left(\frac{-(\overline{H}\theta + \frac{2\gamma T^{(N)}}{\overline{F}})}{2\sqrt{T^{(n)}\gamma}}-\frac{2t\sqrt{T^{(n)}\gamma}}{\theta}\right)\right].
\end{aligned}
\label{eq_mgf_approximated_headway}
\end{equation}
Here $\Phi(\cdot)$ denotes the standard normal cumulative distribution function and its analytic continuation for complex arguments in the root calculation. When $\gamma=0$, no incident duration shocks occur, so $H^{(n)}_{\text{Normal}}=\overline{H}$ almost surely and $M_{H^{(n)}_{\text{Normal}}}(t)=e^{t\overline{H}}$.
\end{proposition}

For compact notation, let $\mu:=\overline{H}+2\gamma T^{(N)}/(\theta\overline{F})$ and $\sigma:=2\sqrt{T^{(n)}\gamma}/\theta$ denote the mean and standard deviation of $\widehat{H}^{(n)}_{\text{Normal}}$. Differentiating the MGF and using $1-\Phi(-\mu/\sigma)=\Phi(\mu/\sigma)$ then gives the corresponding mean and variance of the approximated headway.
\begin{align}
    \mathbb{E}[H^{(n)}_{\text{Normal}}] &= \mu \cdot  \Phi\left(\frac{\mu}{\sigma} \right) +\sigma\cdot \phi\left(\frac{-\mu}{\sigma}\right), \label{eq_mean_normal_tr_hdw}\\
   \operatorname{Var}[H^{(n)}_{\text{Normal}}]&=  \mu \sigma \phi\left(\frac{-\mu}{\sigma}\right) +\Phi\left({\frac{\mu}{\sigma}}\right) \left(\mu^{2} + \sigma^{2}\right) \ - \left( \mu \Phi\left(\frac{\mu}{\sigma}\right) + \phi\left(\frac{-\mu}{\sigma}\right) \sigma\right)^{2},
\end{align}
where $\phi(\cdot)$ is the probability density function of a standard normal distribution. For $\gamma=0$, these expressions are interpreted by continuity, giving mean $\overline{H}$ and variance zero. The following proposition characterizes how incident frequency $\gamma$ and mean incident duration $1/\theta$ affect the expected headway under the normal approximation.

\begin{proposition}[{Impact of incidents on headway}]\label{prop_headway_increase}
Under Assumption \ref{assumption_renewal} and the normal approximation, $\mathbb{E}[H^{(n)}_{\text{Normal}}]$ is nondecreasing in $\gamma$ and $1/\theta$. It is strictly increasing in $\gamma$ and, when $\gamma>0$, strictly increasing in $1/\theta$.
\end{proposition}

Proposition \ref{prop_headway_increase} supports the comparative stability analysis in Section \ref{sec_stable}.

\subsubsection{Distribution of passenger arrivals.}\label{sec_Yn}
We can now return to the passenger arrival distribution that enters the queueing model. Recall that $H^{(n)}_{\text{Normal}}$ is obtained by taking the positive part of the normal approximation to the raw headway $\widehat{H}^{(n)}$. Averaging the conditional Poisson PGF over this headway distribution gives the following result for $Y^{(n)}$.
\begin{proposition}[{PGF of $Y^{(n)}$}]\label{prop_arrival_pax}
For each $n=1,\ldots,N$ and $\gamma>0$, the PGF $Y(z)$ of $Y^{(n)}$ is
\begin{align}
Y(z) = \Phi\left(\frac{-\mu}{\sigma}\right) + e^{\mu \lambda^{(n)}(z-1) + \frac{\sigma^2(\lambda^{(n)}z-\lambda^{(n)})^2}{2}}\left[1 - \Phi\left(\frac{-\mu}{\sigma} - \sigma \lambda^{(n)}(z-1)\right)\right],
\end{align}
where $\mu =   \overline{H} + \frac{2 \gamma T^{(N)}}{\theta\overline{F}}$ and $\sigma = \frac{2 \sqrt{T^{(n)} \gamma}}{\theta}$ are the mean and standard deviation of the raw headway approximation $\widehat{H}^{(n)}_{\text{Normal}}$, respectively.
When $\gamma=0$, $Y^{(n)}$ is Poisson with mean $\lambda^{(n)}\overline{H}$ and $Y(z)=\exp[\lambda^{(n)}\overline{H}(z-1)]$.
\end{proposition}

Differentiating the PGF in Proposition \ref{prop_arrival_pax} at $z=1$ gives the first three factorial moments. We then convert these quantities to $\overline{Y}^{(n)}$, $\overline{\overline{Y}}^{(n)}$, and $\overline{\overline{\overline{Y}}}^{(n)}$. The first moment is shown below. The second and third central moments follow from the same differentiation and conversion, but their longer expressions are omitted for readability.
\begin{align}
\overline{Y}^{(n)} &= \left(\mu \cdot  \Phi\left(\frac{\mu}{\sigma} \right) +\sigma\cdot \phi\left(\frac{-\mu}{\sigma}\right)\right) \cdot \lambda^{(n)}.
\end{align}
When $\gamma=0$, this expression is interpreted by continuity and gives $\overline{Y}^{(n)}=\lambda^{(n)}\overline{H}$.

\subsubsection{Solving for the roots.}\label{sec_root_solving}
With $Y(z)$ now specified, we turn to the $C$ roots of $\mathcal{F}^{(n)}(z)$ in the closed unit disk that are required by the queue length and waiting time formulas. A nonlinear solver with a single initialization generally returns only one root, and repeated initializations may converge to the same location. The expression for $Y(z)$ makes this problem more difficult for the vehicle capacities considered here. We therefore develop a certified interpolation search that combines the structured initialization of \citet{powell1985analysis} with deterministic interpolation and multiplicity certification. The algorithm supplies the required roots and reports a diagnostic failure when certification is unsuccessful.

\textbf{Step 1. Initial root search.} Following the structured angular initialization of \citet{powell1985analysis}, we place $C$ initial values at equally spaced angles $2\pi(k+1/2)/C$ for $k=0,\ldots,C-1$. Let $\rho^{(n)}:=\overline{Y}^{(n)}/\overline{S}^{(n)}$ denote utilization at station $n$. Their common radius is $1-0.5\rho^{(n)}$, bounded between 0.05 and 0.98. From each initial value, a local solver evaluates the real and imaginary parts of $\mathcal{F}^{(n)}(z)$ in Eq. \ref{eq_characteristic_function}. Candidates that pass the convergence, residual, and unit disk checks are retained together with their conjugates and grouped into distinct root clusters. Polar extrapolation from consecutive clusters provides additional initial values for the initial search.

\begin{figure}[htb]
\centering
\includegraphics[width=\linewidth]{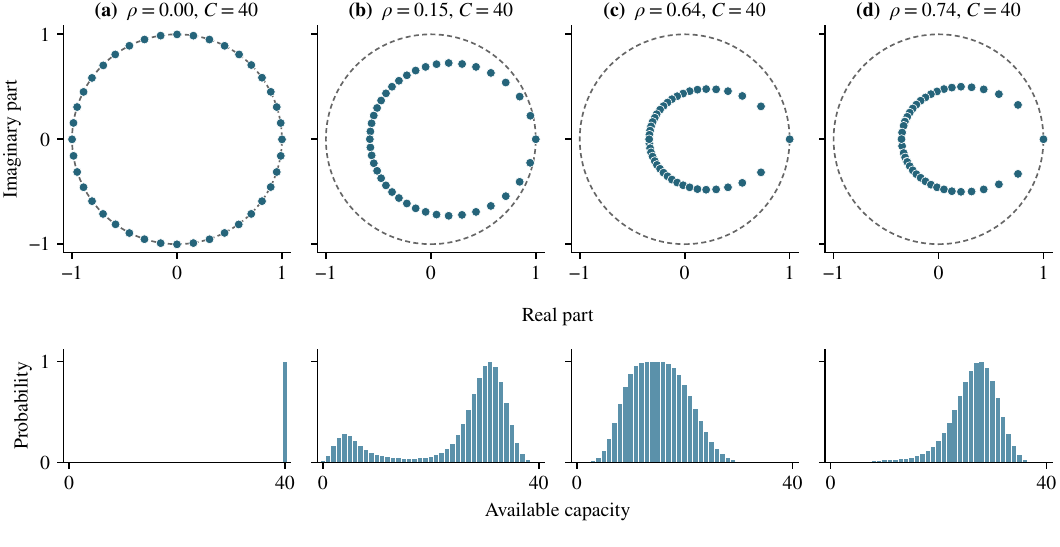}
\caption{Characteristic roots and available capacity distributions across station utilization levels. Each column reports the certified roots and available capacity distribution for the displayed $\rho^{(n)}$ and $C$}
\label{fig_example_root}
\end{figure}

The root locations often form an oval pattern. Figure \ref{fig_example_root} illustrates this pattern for different utilization ratios $\rho^{(n)}$ and available capacity distributions. As $\rho^{(n)}$ increases, the root locations generally move from a nearly circular pattern toward a more elongated pattern.

\begin{figure}[htb]
\centering
\includegraphics[width=\linewidth]{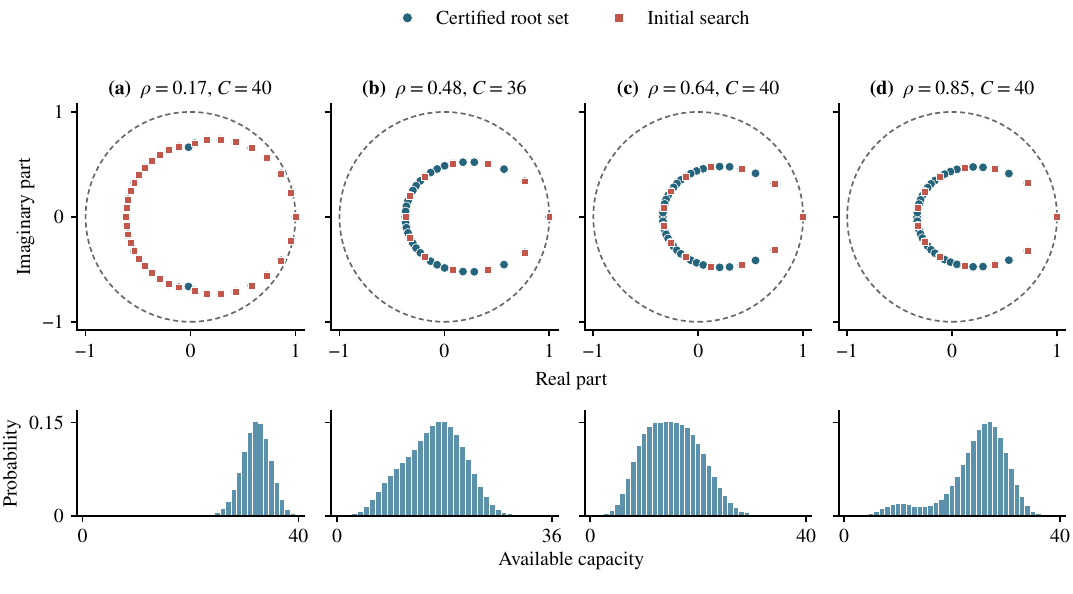}
\caption{Comparison of the initial search and certified root set across station utilization and vehicle capacity settings}
\label{fig_root_step12}
\end{figure}

\textbf{Step 2. Interpolation search.} The initial search may return fewer than $C$ roots counted according to multiplicity. We order the distinct root clusters by complex argument and interpret adjacent pairs circularly. At interpolation level $L$, the algorithm interpolates both the modulus and the unwrapped argument of every adjacent pair at $L-1$ locations. It places initial values at the interpolated modulus and at radial factors 0.85 and 1.15, with the modulus bounded by one. These radial values prevent a straight chord between two roots from moving all initial values toward the origin. The algorithm applies the local solver, refines each successful candidate, and updates the root clusters. If an iteration finds no new root location, it increases $L$. Because the coefficients of $\mathcal{F}^{(n)}(z)$ are real, the conjugate of every nonreal candidate is also evaluated. Figure \ref{fig_root_step12} compares the distinct locations from the initial search with the root set after interpolation and certification.

\textbf{Step 3. Multiplicity certification and termination.} Finding $C$ distinct locations is not a valid stopping rule because a root may have multiplicity greater than one. Let $\mathcal{Z}^{(k)} := \{\zeta_0^{(k)},\ldots,\zeta_{M_k-1}^{(k)}\}$ contain the $M_k$ distinct root clusters after iteration $k$. For cluster $j$, define the counterclockwise circular contour $\Gamma_j^{(k)} := \{\zeta_j^{(k)}+r_j^{(k)}\exp(\varphi\sqrt{-1})\mid0\leq\varphi\leq2\pi\}$. Its radius is smaller than half the distance to the nearest other cluster and no greater than 0.05. We use 0.05 when only one cluster is available. Thus, $\Gamma_j^{(k)}$ is tied to candidate $\zeta_j^{(k)}$ and separates it from other candidates. We use the contour only when $\mathcal{F}^{(n)}(z)\neq0$ for every $z$ on $\Gamma_j^{(k)}$.
\enlargethispage{\baselineskip}
\begin{align}
m_j^{(k)} := \frac{1}{2\pi\sqrt{-1}}\int_{\Gamma_j^{(k)}}\frac{[\mathcal{F}^{(n)}]'(z)}{\mathcal{F}^{(n)}(z)}\,dz.
\label{eq_root_multiplicity}
\end{align}
Because $\mathcal{F}^{(n)}(z)$ has no poles, $m_j^{(k)}$ equals the number of roots enclosed by $\Gamma_j^{(k)}$, counted according to multiplicity. It is also the winding number of $\mathcal{F}^{(n)}(z)$ around the origin along this contour. A value of one certifies a simple root in cluster $j$, while a value of two certifies a root of multiplicity two. We compute this value from the unwrapped change in the argument of $\mathcal{F}^{(n)}(z)$ and repeat the calculation with finer discretizations and smaller contours. A cluster is certified when these calculations return the same positive integer.

The algorithm terminates when every cluster is certified and the sum of the certified multiplicities equals $C$. Otherwise, it returns to Step 2 until a search limit is reached and then reports a diagnostic failure. This certification follows numerical treatments of transform equations with possible multiple roots \citep{oblakova2019exact}.

The complete procedure is provided in Algorithm \ref{alg_inter_search} of the EC appendix.

\subsection{Stability condition}\label{sec_stable}
The preceding transform calculations characterize the queue conditional on the existence of a steady state. We now establish when this steady state exists directly from the transient queue recursion. We call a station stable when its queue process admits a proper stationary probability distribution.

Assumption \ref{assumption_station_service} separates the propagation of available capacity along the route from the evolution of the queue at an individual station. Upstream passenger flows determine $s^{(n)}$, so available capacity remains endogenous across stations. After $s^{(n)}$ is propagated to station $n$, the model uses its complete distribution to generate repeated service opportunities that are exogenous to the local queue, thereby carrying the effect of upstream crowding into the stationary single station analysis.

Assumption \ref{assumption_station_service} is exact at the first station because every arriving vehicle has capacity $C$. Under persistent upstream saturation, the approximation is also reasonable because consecutive vehicles depart at capacity and independent binomial alighting determines the available capacity at station $n$. Outside these cases, independent renewal headways generate new passenger arrivals for successive vehicles, while alighting at intermediate stations weakens the persistence of upstream loads. These mechanisms limit dependence in available capacity when suspensions are short and rarely affect several vehicle runs. When long disruptions or common operating conditions sustain correlated vehicle loads, the approximation becomes weaker. Section \ref{sec_num_ex} assesses it through recursive simulation.

\begin{proposition}[{Stability condition}]\label{prop_stability}
Under Assumptions \ref{assumption_renewal} and \ref{assumption_station_service} and the normal approximation in Proposition \ref{prop_headway_distribution}, the bulk service queueing system at station $n$ with $\overline{S}^{(n)}>0$ is stable if and only if
\begin{align}
\rho^{(n)} = \frac{\overline{Y}^{(n)}}{\overline{S}^{(n)}} = \frac{\left(\mu \cdot  \Phi\left(\frac{\mu}{\sigma} \right) +\sigma\cdot \phi\left(\frac{-\mu}{\sigma}\right)\right) \cdot \lambda^{(n)} }{\sum_{u = 0}^C s_u^{(n)}u} = \frac{\lambda^{(n)} \cdot \mathbb{E}[H^{(n)}_{\text{Normal}}]}{\sum_{u = 0}^C s_u^{(n)}u}< 1,
\label{eq_stab}
\end{align}
where $\rho^{(n)}$ is the utilization ratio at station $n$.
When $\gamma=0$, the expression is interpreted by continuity, giving $\mathbb{E}[H^{(n)}_{\text{Normal}}]=\overline{H}$ and $\rho^{(n)}=\lambda^{(n)}\overline{H}/\sum_{u=0}^C s_u^{(n)}u$.
\end{proposition}

The proof in the EC appendix writes the residual queue as a Lindley recursion \citep{lindley1952theory} and applies the stability result of \citet{loynes1962stability}. It therefore establishes the existence of the stationary distribution without using $q_0^{(n)}$ or any other quantity derived under stationarity.

Proposition \ref{prop_stability} states that station $n$ is stable when the mean number of passenger arrivals within a headway is smaller than the mean available capacity of an arriving vehicle after alighting. Proposition \ref{prop_headway_increase} shows that a higher incident rate or a longer mean incident duration increases $\mathbb{E}[H^{(n)}_{\text{Normal}}]$ under the conditions stated there, making instability more likely. The result therefore quantifies the throughput loss induced by incidents.

\begin{remark}
The distribution $s^{(n)}$ depends on the upstream queue and vehicle load distributions, which in turn require the characteristic roots at stations $1,\ldots,n-1$. Stability at station $n$ must therefore be evaluated after propagating the preceding stations. At station 1, $s_C^{(1)}=1$ and $s_u^{(1)}=0$ for $u=0,\ldots,C-1$, so Eq. \ref{eq_stab} gives $\rho^{(1)}=\lambda^{(1)}\mathbb{E}[H^{(1)}_{\text{Normal}}]/C$ for a direct stability assessment.
\end{remark}

\begin{remark}
Proposition \ref{prop_stability} concerns station level stability. A route is stable if and only if $\rho^{(n)}<1$ for every $n=1,\ldots,N$.
\end{remark}

\begin{remark}
If station $n-1$ is stable, its departing load distribution determines $s^{(n)}$ through Section \ref{sec_space_dist}, after which stability at station $n$ can be evaluated. If station $n-1$ is unstable, station $n$ may nevertheless be stable because passengers can alight there. In this case, vehicles depart station $n-1$ at capacity, so $v_C^{(n-1)}=1$ and $v_k^{(n-1)}=0$ for $k=0,\ldots,C-1$. Binomial alighting then gives $\overline{S}^{(n)}=\alpha^{(n)}C$, and the stability condition becomes $\rho^{(n)}=\lambda^{(n)}\mathbb{E}[H^{(n)}_{\text{Normal}}]/[\alpha^{(n)}C]<1$.
\end{remark}

\subsection{Summary of calculation procedure}\label{sec_summary_calculation}
Algorithm \ref{alg_overall} combines the preceding analytical results into a sequential calculation over the $N$ stations of the route.

\begin{algorithm}[htb]
\caption{Calculation of performance indicators} \label{alg_overall}

\fontsize{9.5}{11.4}\selectfont
\begin{algorithmic}[1]
\State Initialize $v^{(0)}_0\gets1$ and $v^{(0)}_k\gets0$ for $k=1,\ldots,C$
\For{$n= 1:N$}
    \State $g^{(n)} = v^{(n-1)}A^{(n)}$  \Comment{Proposition \ref{prop_space_dist}}
    \State $s_k^{(n)}\gets g_{C-k}^{(n)}$ for $k=0,1,\ldots,C$ \Comment{Eq. \ref{eq_s}}
    \State Calculate $\overline{S}^{(n)}, \overline{\overline{S}}^{(n)},$ and $\overline{\overline{\overline{S}}}^{(n)}$ based on $s^{(n)}$. 
    \State Calculate $\overline{Y}^{(n)}, \overline{\overline{Y}}^{(n)},$ and $\overline{\overline{\overline{Y}}}^{(n)}$  \Comment{Section \ref{sec_Yn}}
    \If{$\overline{Y}^{(n)} <  \overline{S}^{(n)}$}  \Comment{Station $n$ is stable}
    \State Solve for the roots $z_0^*,\ldots,z^*_{C-1}$ of the denominator of $Q(z)$ in Eq. \ref{eq_Qz1}\Comment{Section \ref{sec_root_solving}}
    \State Calculate $q_0^{(n)},\ldots,q_{C-1}^{(n)}$ from $z_0^*,\ldots,z^*_{C-1}$ \Comment{Section \ref{sec_cal_q}}
    \State Calculate $\mathbb{E}[Q^{(n)}]$ and $\operatorname{Var}[Q^{(n)}]$, and, if $\lambda^{(n)}>0$, calculate $\mathbb{E}[W^{(n)}]$ and $\operatorname{Var}[W^{(n)}]$ \Comment{Eqs. \ref{eq_EQ}--\ref{eq_VarW}}
    \State $v^{(n)} = g^{(n)}B^{(n)}$ \Comment{Eq. \ref{eq_v}. $B^{(n)}$ is a function of $q_{0:C-1}^{(n)}$}
    \Else \Comment{Station $n$ is not stable}
    \State Set $\mathbb{E}[Q^{(n)}], \operatorname{Var}[Q^{(n)}], \mathbb{E}[W^{(n)}]$, and $\operatorname{Var}[W^{(n)}]$ to infinity
    \State Set $v^{(n)}_C\gets1$ and $v^{(n)}_k\gets0$ for $k=0,1,\ldots,C-1$
    \EndIf
\EndFor
\end{algorithmic}
\end{algorithm}

\section{Numerical experiments}\label{sec_num_ex}

This section illustrates the analytical framework on a bus route with ten stations and assesses the resulting resilience indicators against recursive FIFO simulation. The experiments isolate the effects of vehicle capacity, incident frequency, incident duration, scheduled headway, and demand, with particular attention to stations with limited capacity slack. All analytical results use the certified interpolation root search described in Section \ref{sec_root_solving}. Section \ref{append_root_search} of the EC appendix reports the numerical settings used for the search.

\subsection{Experimental design}\label{sec_experiment}
\begin{figure}[htb]
\centering
\subfloat{\includegraphics[width=0.75\textwidth]{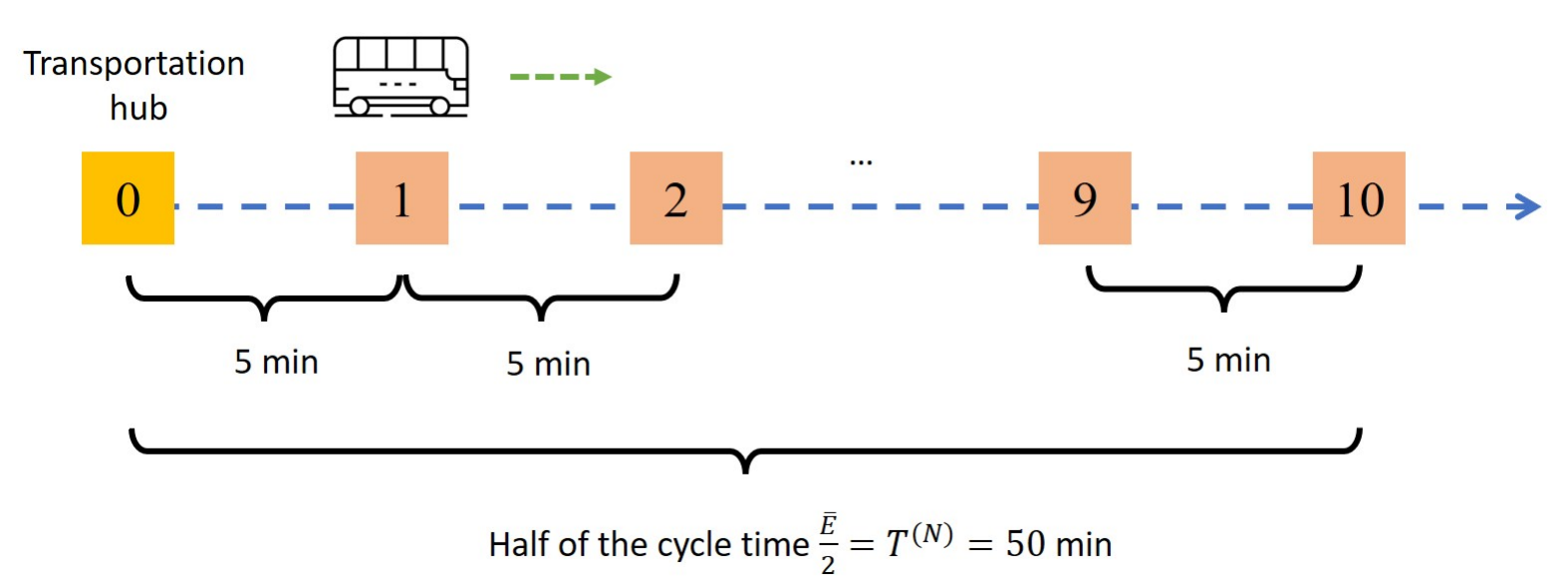}}
\caption{Case study route layout}
\label{fig_case_layout}
\end{figure}

We use an example bus route adapted from \citet{islam2015model}. Table \ref{tab_bus_route} reports the attributes of its ten stations. Figure \ref{fig_case_layout} shows the route layout. Travel time without incidents is 5 minutes between consecutive stations, total cycle time is $\overline{E}=100$ minutes, and travel time from the hub to the last station is $T^{(N)}=50$ minutes.
\begin{table}[htb]
\centering
\caption{Example bus system parameters.  $\lambda^{(n)}$ and $\alpha^{(n)}$ are the passenger arrival rate and alighting probability at station $n$}
\fontsize{8.6}{10.3}\selectfont
\setlength{\tabcolsep}{5pt}
\begin{tabular}{@{}cccccc@{}}
\toprule
Station & \makecell{$\lambda^{(n)}$\\(passengers/min)} & $\alpha^{(n)}$ & Station & \makecell{$\lambda^{(n)}$\\(passengers/min)} & $\alpha^{(n)}$ \\ \midrule
1          & 0.75                      & 0              & 6          & 1                         & 0.8            \\
2          & 2                         & 0              & 7          & 0.75                      & 0.5            \\
3          & 0.75                      & 0.1            & 8          & 1.5                       & 0.3            \\
4          & 0.5                       & 0.25           & 9          & 0.2                       & 0.75           \\
5          & 0.5                       & 0.25           & 10         & 0                         & 1              \\ \bottomrule
\end{tabular}
\label{tab_bus_route}
\end{table}

To test the sensitivity of resilience indicators to different parameters, we consider different values of $C$, $\theta$, $\gamma$, $\overline{H}$, and demand in Table \ref{tab_scenario}. The demand factor multiplies the base arrival rates $\lambda^{(n)}$ in Table \ref{tab_bus_route}. Keeping the route cycle time at $\overline{E}=100$ minutes, we select the integer fleet sizes $\overline{F}\in\{50,25,14\}$, which yield $\overline{H}=\overline{E}/\overline{F}\in\{2,4,50/7\}$ minutes. The reference scenario sets $\overline{F}=25$ and $\overline{H}=4$ minutes; each test varies one parameter and retains all other reference values.

\begin{table}[htb]
\centering
\caption{Scenario design for vehicle capacity $C$, incident rate $\gamma$, incident duration rate $\theta$, scheduled headway $\overline{H}$, and demand factor}
\fontsize{8.6}{10.3}\selectfont
\setlength{\tabcolsep}{7pt}
\begin{tabular}{@{}ccc@{}}
\toprule
Parameters      & Tested values                             & Reference value \\ \midrule
$C$             & \{30, 34, 38\}                      & 34              \\
$\gamma$ (/min)        & \{0, 1/10, 1/5, 1/3\} & 1/5                \\
$\theta$ (/min)       & \{2, 1 ,1/2\} & 1              \\
$\overline{H}$ (min) & \{2, 4, 50/7\}                             & 4               \\
Demand factor   & \{0.5, 0.75, 1\}                & 0.75             \\ \bottomrule
\end{tabular}
\label{tab_scenario}
\end{table}
\subsection{Performance indicators}

\begin{figure}[htb]
\captionsetup[subfigure]{justification=centering}
\centering
\subfloat[Sensitivity on $C$]{\includegraphics[width=0.33\linewidth]{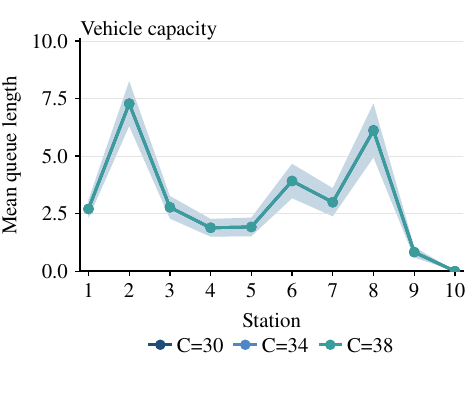} \label{fig_sen_C_on_Q}}
\subfloat[Sensitivity on $\gamma$]{\includegraphics[width=0.33\linewidth]{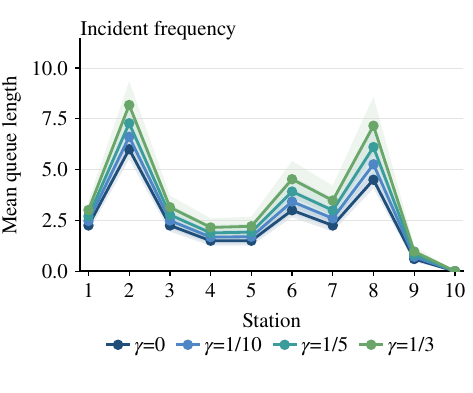}\label{fig_sen_gamma_on_Q}}
\subfloat[Sensitivity on $\theta$]{\includegraphics[width=0.33\linewidth]{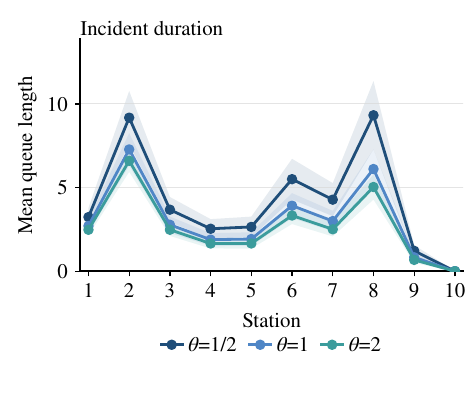}\label{fig_sen_theta_on_Q}}\\
\subfloat[Sensitivity on $\overline{H}$]{\includegraphics[width=0.33\linewidth]{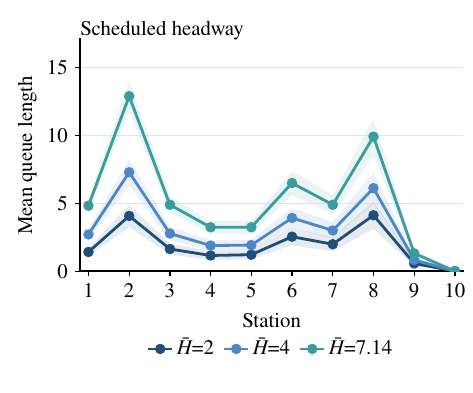}\label{fig_sen_H_bar_on_Q}}
\subfloat[Sensitivity on demand factor]{\includegraphics[width=0.33\linewidth]{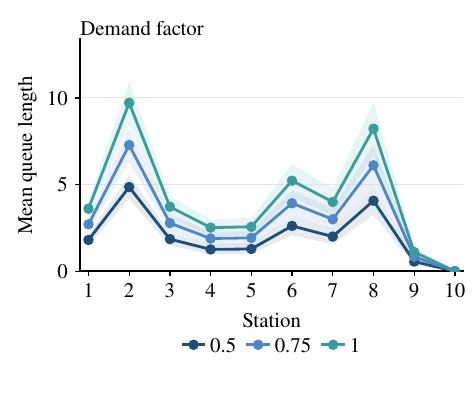}\label{fig_sen_demand_on_Q}}
\caption{Mean queue length across stations. Lines show the means, and shaded bands extend $0.2$ standard deviations above and below each mean}
\label{fig_change_on_Q}
\end{figure}

Figure \ref{fig_change_on_Q} shows the mean and standard deviation of queue length under the tested scenarios. The patterns are consistent with the passenger arrival and alighting rates. Expected queues are relatively high at stations 2 and 8, whereas the expected queue at station 10 is zero because $\lambda^{(10)}=0$.

Figure \ref{fig_sen_C_on_Q} shows limited sensitivity to vehicle capacity within the tested range because capacity is not fully utilized under the reference setting. Figure \ref{fig_sen_gamma_on_Q} shows the effect of incident occurrence rate $\gamma$. When $\gamma=0$, the expected queue length at station 8 is 4.5 passengers. As $\gamma$ increases, expected queue lengths and their standard deviations rise. At $\gamma=1/3$ per minute, the expected queue at station 8 reaches 7.2 passengers. Incident duration produces a similar pattern in Figure \ref{fig_sen_theta_on_Q}. When mean incident duration is 30 seconds and $\theta=2$, $\mathbb{E}[Q^{(8)}]=5.0$. When mean incident duration is 2 minutes and $\theta=1/2$, this expectation increases to 9.3. Both effects are larger at crowded stations.

Figure \ref{fig_sen_H_bar_on_Q} shows that a longer scheduled headway lowers the service rate and increases the expected queue. As $\overline{H}$ increases from 2 minutes to $50/7$ minutes, the expected queue at station 8 increases from 4.1 to 9.9 passengers. Figure \ref{fig_sen_demand_on_Q} shows a similar pattern for demand. Increasing the demand factor from 0.5 to 1.0 raises the expected queue at station 8 from 4.1 to 8.2 passengers. The effects of $\overline{H}$ and demand are more uniform across stations than the effects of incident frequency and duration.

\begin{figure}[htb]
\captionsetup[subfigure]{justification=centering}
\centering
\subfloat[Sensitivity on $C$]{\includegraphics[width=0.33\linewidth]{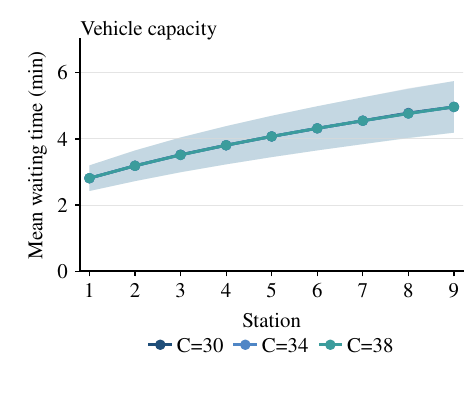} \label{fig_sen_C_on_W}}
\subfloat[Sensitivity on $\gamma$]{\includegraphics[width=0.33\linewidth]{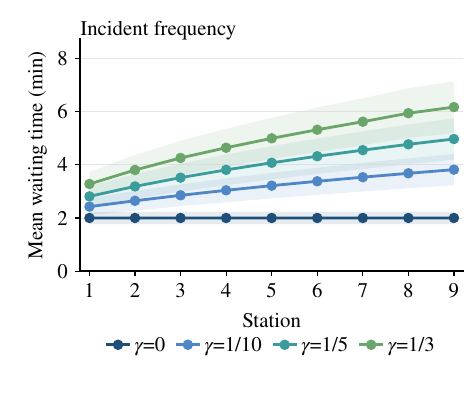}\label{fig_sen_gamma_on_W}}
\subfloat[Sensitivity on $\theta$]{\includegraphics[width=0.33\linewidth]{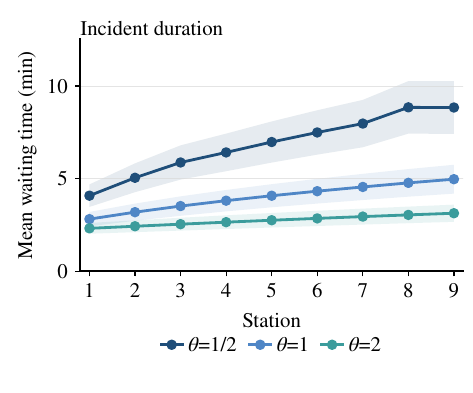}\label{fig_sen_theta_on_W}}\\
\subfloat[Sensitivity on $\overline{H}$]{\includegraphics[width=0.33\linewidth]{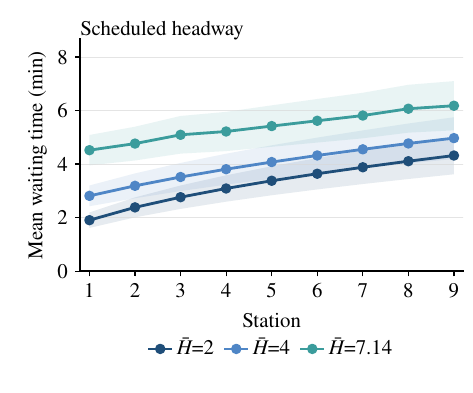}}
\subfloat[Sensitivity on demand factor]{\includegraphics[width=0.33\linewidth]{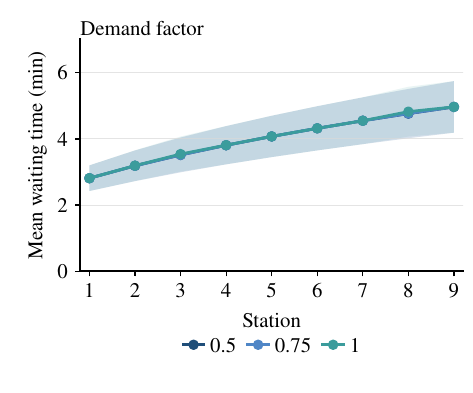}}
\caption{Mean waiting time across stations. Lines show the means, and shaded bands extend $0.2$ standard deviations above and below each mean}
\label{fig_change_on_W}
\end{figure}

Figure \ref{fig_change_on_W} shows the mean and standard deviation of passenger waiting time. Station 10 is omitted because its zero passenger arrival rate produces no waiting time observations. Downstream stations generally have higher waiting time means and standard deviations because headway variability accumulates along the route. At congested stations such as stations 3 and 8, capacity constraints leave some passengers behind and further increase waiting time.

Figure \ref{fig_sen_C_on_W} shows that capacity has a limited effect on waiting time within the tested range. Figures \ref{fig_sen_gamma_on_W} and \ref{fig_sen_theta_on_W} show the effects of $\gamma$ and $\theta$. Higher values of $\gamma$ and $1/\theta$ increase the expected headway and, in turn, mean waiting time at stations with passenger demand, particularly at crowded stations. When $\gamma=0$, there are no incidents, left behind passengers, or headway irregularities. Expected waiting time at stations 1 through 9 is then $\overline{H}/2=2$ minutes under the reference headway of 4 minutes. When $\gamma$ increases to $1/5$, station 3 has left behind passengers and its expected waiting time increases to 3.5 minutes. When $\theta$ decreases from 2 to $1/2$, expected waiting time at station 8 increases from 3.0 to 8.8 minutes.

A larger scheduled headway directly increases waiting time at stations with passenger demand, and passengers are left behind at stations 3 and 8 when $\overline{H}=50/7$ minutes. Demand does not alter the headway distribution, so its effect on waiting time remains limited until capacity constraints leave passengers behind. At station 8, increasing the demand factor from 0.5 to 1.0 raises expected waiting time only from 4.77 to 4.83 minutes within the tested range.

\subsection{Comparison between simulation and analytical results}

To assess the renewal headway, marginal route decoupling, and normal distribution approximations, we develop a recursive FIFO simulation model for the expectation and variance of queue length and waiting time. Unlike the analytical renewal model, the simulation carries the effect of a delayed vehicle to its successors through realized departure times and jointly evolves passenger queues and vehicle loads; see Section \ref{append_sim_alg} of the EC appendix.

\begin{figure}[htb]
\captionsetup[subfigure]{justification=centering}
\centering
\subfloat[\text{$\mathbb{E}[Q^{(n)}]$}]{\includegraphics[width=0.47\linewidth]{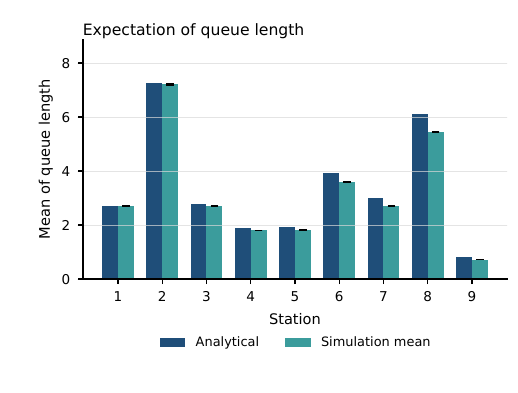}}
\subfloat[Standard deviation of $Q^{(n)}$]{\includegraphics[width=0.47\linewidth]{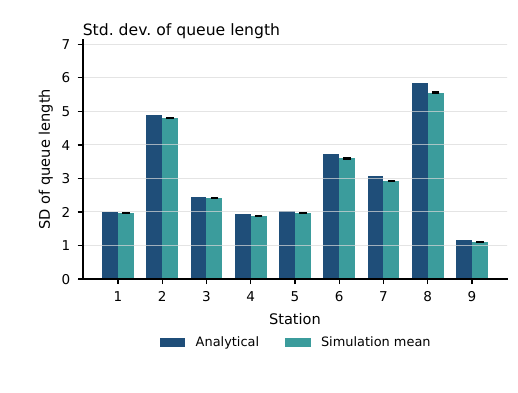}}\\[-2pt]
\subfloat[\text{$\mathbb{E}[W^{(n)}]$}]{\includegraphics[width=0.47\linewidth]{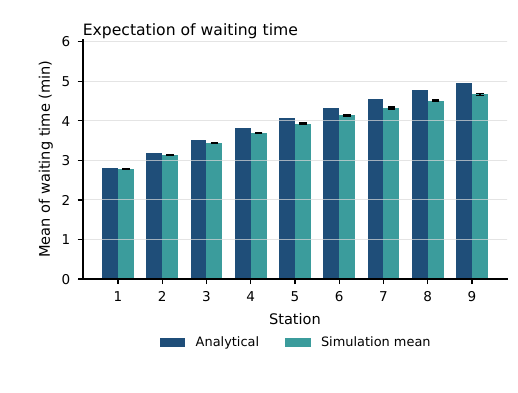}}
\subfloat[Standard deviation of $W^{(n)}$]{\includegraphics[width=0.47\linewidth]{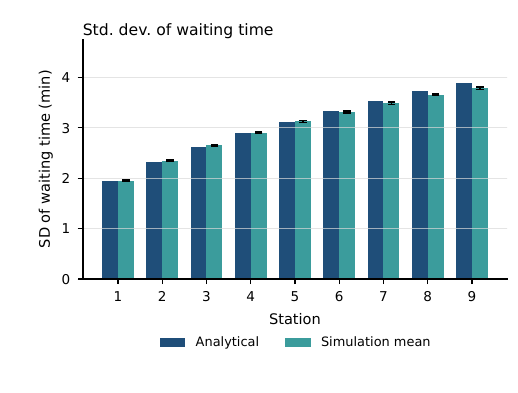}}
\caption{Comparison between simulation and analytical results for the reference setting.}
\label{fig_comp_to_sim}
\end{figure}

The assessment covers seven settings from Table \ref{tab_scenario}, including the reference setting and variations in capacity, headway, demand, incident frequency, and incident duration. For each setting, we conduct 20 independent replications of 10,000 vehicle runs using distinct random seeds and exclude the first 10\% of vehicle runs. Figure \ref{fig_comp_to_sim} compares the analytical results with the means across replications for the reference setting. The whiskers give 95\% confidence intervals (CIs) based on variation across replications. Station 10 has $\lambda^{(10)}=0$, so no passenger arrives there and no waiting time observation exists. We therefore report waiting time through station 9. All error calculations also use stations 1 through 9 so that the four performance measures have common support.

Table \ref{tab_sim_validation} reports the relative $\ell_1$ error in percentage terms and the average half width (HW) of the station CIs. For a given performance metric, $x_{\mathrm{Sim}}^{(n)}$ denotes the mean estimate across the 20 independent simulation replications at station $n$, and $x_{\mathrm{Ana}}^{(n)}$ denotes the corresponding result obtained from the analytical model. The relative $\ell_1$ error is calculated as $100\sum_{n=1}^{9}|x_{\mathrm{Sim}}^{(n)}-x_{\mathrm{Ana}}^{(n)}|/\sum_{n=1}^{9}|x_{\mathrm{Ana}}^{(n)}|$. Under the reference setting, the queue length errors are 5.64\% for the mean and 3.27\% for the standard deviation. The corresponding waiting time errors are 3.87\% and 1.26\%.

Across all tested settings, the waiting time errors remain below 8.21\%. The largest discrepancy is the 16.58\% error in mean queue length when the mean incident duration is 2 minutes. This setting is also where the short suspension premise and the independent renewal approximation are most strained. The average confidence interval half width over the seven settings is at most 0.013 passengers for mean queue length, 0.025 passengers for its standard deviation, 0.039 minutes for mean waiting time, and 0.052 minutes for its standard deviation. The sampling uncertainty is therefore small relative to the main analytical differences.

\begin{table}[htb]
\centering
\begin{threeparttable}
\caption{Relative $\ell_1$ errors and confidence interval half widths for queue length and waiting time}
\vspace{-4pt}
\fontsize{8.6}{10.3}\selectfont
\begin{tabular*}{\linewidth}{@{\extracolsep{\fill}}lcccccccc@{}}
\toprule
& \multicolumn{2}{c}{$\mathbb{E}[Q^{(n)}]$} & \multicolumn{2}{c}{$\mathrm{SD}[Q^{(n)}]$} & \multicolumn{2}{c}{$\mathbb{E}[W^{(n)}]$} & \multicolumn{2}{c}{$\mathrm{SD}[W^{(n)}]$} \\
\cmidrule(lr){2-3}\cmidrule(lr){4-5}\cmidrule(lr){6-7}\cmidrule(lr){8-9}
Setting & Error (\%) & CI HW & Error (\%) & CI HW & Error (\%) & CI HW & Error (\%) & CI HW \\ \midrule
$C=30$ & 5.40 & 0.008 & 3.34 & 0.012 & 3.88 & 0.016 & 1.32 & 0.020 \\
Reference & 5.64 & 0.009 & 3.27 & 0.013 & 3.87 & 0.016 & 1.26 & 0.018 \\
$\overline{H}=50/7$ & 0.67 & 0.012 & 1.49 & 0.015 & 0.93 & 0.013 & 0.31 & 0.016 \\
Demand factor $=1$ & 4.89 & 0.010 & 3.53 & 0.016 & 4.04 & 0.014 & 1.38 & 0.016 \\
$\gamma=0$ & 0.11 & 0.008 & 0.14 & 0.006 & 0.05 & 0.004 & 0.07 & 0.002 \\
$\gamma=1/3$ & 7.97 & 0.010 & 4.85 & 0.016 & 5.48 & 0.020 & 2.68 & 0.024 \\
$\theta=1/2$ & 16.58 & 0.011 & 8.68 & 0.025 & 8.20 & 0.039 & 3.65 & 0.052 \\ \bottomrule
\end{tabular*}
\label{tab_sim_validation}
\begin{tablenotes}[flushleft]
\fontsize{8.2}{9.8}\selectfont
\item[] The confidence interval half widths are averaged across stations 1 through 9. Queue length is measured in passengers, and waiting time is measured in minutes. In each nonreference row, only the listed parameter differs from its reference value.
\end{tablenotes}
\end{threeparttable}
\vspace{-4pt}
\end{table}

The results provide numerical evidence for the accuracy of the approximations within the tested parameter range, but they do not establish exact equivalence between the analytical renewal process and the recursive FIFO process. One source of discrepancy is the normal approximation of the raw headway distribution. As shown in Figure \ref{fig_emp_valid}, the compound Poisson based raw headway concentrates more probability density near the mean than the normal approximation and therefore assigns a lower probability to large deviations from the scheduled headway. A second source is the independent redraw of each analytical headway, whereas the recursive FIFO simulation carries delay forward.

\section{Conclusion and discussion}\label{sec_conclusion}
This paper proposes a stochastic framework for quantifying the resilience of PTSs under short random service suspensions. Specifically, we introduce an independent renewal approximation for nonnegative headways and derive stability conditions and closed form expressions for the mean and variance of queue length and waiting time at each station. The common headway marginal preserves the incident effects implied by a two state vehicle speed process, while the renewal and route decoupling assumptions make the bulk service analysis tractable. The stability conditions imply that the system is more likely to be unstable under high incident rates, long incident durations, high demand, low service frequency, and low vehicle capacity. We implement the model on an example bus route adapted from the literature, and the sensitivity analysis shows that congested stations are more vulnerable to random service suspensions. Comparison with a recursive FIFO simulation across seven parameter settings quantifies the approximation error and shows that the discrepancy increases when incident durations move beyond the intended short suspension setting.

The proposed model can support public transit planning, operations, and management in several ways. For the design and planning of resilient systems under random system interruptions, the model can inform headway design and vehicle capacity determination. Its queue length estimates can be used to evaluate the layout and capacity of congested stations. Historical demand and incident information can be used with the model to monitor system resilience and identify critical stations. The estimates of waiting time and queue length can also support cost benefit analyses of service improvements. For example, the model can identify whether a larger vehicle, shorter headway, or higher maintenance frequency is the most cost effective way to keep waiting time below a threshold.

Several extensions are possible. A network level model would differ primarily through the inclusion of transfer passengers. Transfer demand could enter the arrival process, but the model would require additional parameters related to transfers and alighting rates. The renewal approximation could also be replaced by a dependence aware model that carries recursive FIFO delay and the joint passenger queue and vehicle load state. Another extension would consider the partial interruptions studied in Section \ref{sec_liter_disrupt}. Vehicles retain a positive speed during a partial interruption, which requires a revised headway distribution. A final extension would relax the assumption of no passenger balking or reneging and incorporate more complicated passenger behavior.

\ACKNOWLEDGMENT{This work was in part supported by SJTU UM Joint Institute, J. Wu \& J. Sun Endowment Fund, C2SMART University Transportation Center, and NYU Tandon School of Engineering.}

\def\bibsep{0pt}
\begingroup
\exhyphenpenalty=10000
\bibliography{mybibfile}

@article{rahimi2019analysis,
  title={Analysis of transit users’ waiting tolerance in response to unplanned service disruptions},
  author={Rahimi, Ehsan and Shamshiripour, Ali and Shabanpour, Ramin and Mohammadian, Abolfazl and Auld, Joshua},
  journal={Transportation Research Part D: Transport and Environment},
  volume={77},
  pages={639--653},
  year={2019},
  publisher={Elsevier}
}

@article{mo2020capacity,
  title={Capacity-constrained network performance model for urban rail systems},
  author={Mo, Baichuan and Ma, Zhenliang and Koutsopoulos, Haris N and Zhao, Jinhua},
  journal={Transportation Research Record},
  volume={2674},
  number={5},
  pages={59--69},
  year={2020},
  publisher={SAGE Publications Sage CA: Los Angeles, CA}
}

@article{krishnamoorthy2014queues,
  title={Queues with interruptions: a survey},
  author={Krishnamoorthy, Achyutha and Pramod, Padinhare K and Chakravarthy, Srinivas R},
  journal={Top},
  volume={22},
  number={1},
  pages={290--320},
  year={2014},
  publisher={Springer}
}

@article{bailey1954queueing,
  title={On queueing processes with bulk service},
  author={Bailey, Norman TJ},
  journal={Journal of the Royal Statistical Society: Series B (Methodological)},
  volume={16},
  number={1},
  pages={80--87},
  year={1954},
  publisher={Wiley Online Library}
}

@article{madan1989single,
  title={A single channel queue with bulk service subject to interruptions},
  author={Madan, KC},
  journal={Microelectronics Reliability},
  volume={29},
  number={5},
  pages={813--818},
  year={1989},
  publisher={Elsevier}
}

@article{islam2014bulk,
  title={A bulk queue model for the evaluation of impact of headway variations and passenger waiting behavior on public transit performance},
  author={Islam, Md Kamrul and Vandebona, Upali and Dixit, Vinayak V and Sharma, Ashish},
  journal={IEEE Transactions on Intelligent Transportation Systems},
  volume={15},
  number={6},
  pages={2432--2442},
  year={2014},
  publisher={IEEE}
}

@article{mo2022impact,
  title={Impact of unplanned long-term service disruptions on urban public transit systems},
  author={Mo, Baichuan and Von Franque, Max Y and Koutsopoulos, Haris N and Attanucci, John P and Zhao, Jinhua},
  journal={IEEE Open Journal of Intelligent Transportation Systems},
  volume={3},
  pages={551--569},
  year={2022},
  publisher={IEEE}
}

@article{mo2023robust,
  title={Robust path recommendations during public transit disruptions under demand uncertainty},
  author={Mo, Baichuan and Koutsopoulos, Haris N and Shen, Zuo-Jun Max and Zhao, Jinhua},
  journal={Transportation Research Part B: Methodological},
  volume={169},
  pages={82--107},
  year={2023},
  publisher={Elsevier}
}

@article{mo2025individual,
  title={Individual path recommendation under public transit service disruptions considering behavior uncertainty},
  author={Mo, Baichuan and Koutsopoulos, Haris N and Shen, Zuo-Jun Max and Zhao, Jinhua},
  journal={Transportation Science},
  volume={59},
  number={6},
  pages={1235--1258},
  year={2025},
  publisher={INFORMS}
}

@phdthesis{powell1981stochastic,
  title={Stochastic delays in transportation terminals: New results in the theory and application of bulk queues},
  author={Powell, Warren Buckler},
  year={1981},
  school={Massachusetts Institute of Technology}
}

@article{powell1985analysis,
  title={Analysis of vehicle holding and cancellation strategies in bulk arrival, bulk service queues},
  author={Powell, Warren B},
  journal={Transportation Science},
  volume={19},
  number={4},
  pages={352--377},
  year={1985},
  publisher={INFORMS}
}

@article{adan2006rouche,
  title={On the application of {Rouch\'{e}}'s theorem in queueing theory},
  author={Adan, Ivo J B F and van Leeuwaarden, Johan S H and Winands, Emiel M M},
  journal={Operations Research Letters},
  volume={34},
  number={3},
  pages={355--360},
  year={2006},
  publisher={Elsevier}
}

@article{oblakova2019exact,
  title={An exact root-free method for the expected queue length for a class of discrete-time queueing systems},
  author={Oblakova, Anna and Al Hanbali, Ahmad and Boucherie, Richard J and van Ommeren, Jan-Kees W and Zijm, W H M},
  journal={Queueing Systems},
  volume={92},
  number={3--4},
  pages={257--292},
  year={2019},
  doi={10.1007/s11134-019-09614-1},
  publisher={Springer}
}

@article{chaudhry1987computational,
  title={Computational analysis of steady-state probabilities of M/G a, b/1 and related nonbulk queues},
  author={Chaudhry, Mohan L and Madill, BR and Briere, G},
  journal={Queueing systems},
  volume={2},
  number={2},
  pages={93--114},
  year={1987},
  publisher={Springer}
}

@article{bellei2010transit,
  title={Transit vehicles’ headway distribution and service irregularity},
  author={Bellei, Giuseppe and Gkoumas, Konstantinos},
  journal={Public transport},
  volume={2},
  number={4},
  pages={269--289},
  year={2010},
  publisher={Springer}
}

@article{daganzo2009headway,
  title={A headway-based approach to eliminate bus bunching: Systematic analysis and comparisons},
  author={Daganzo, Carlos F},
  journal={Transportation Research Part B: Methodological},
  volume={43},
  number={10},
  pages={913--921},
  year={2009},
  publisher={Elsevier}
}

@article{hickman2001analytic,
  title={An analytic stochastic model for the transit vehicle holding problem},
  author={Hickman, Mark D},
  journal={Transportation Science},
  volume={35},
  number={3},
  pages={215--237},
  year={2001},
  publisher={INFORMS}
}

@article{jin2016optimizing,
  title={Optimizing bus bridging services in response to disruptions of urban transit rail networks},
  author={Jin, Jian Gang and Teo, Kwong Meng and Odoni, Amedeo R},
  journal={Transportation Science},
  volume={50},
  number={3},
  pages={790--804},
  year={2016},
  publisher={INFORMS}
}

@article{osuna1972control,
  title={Control strategies for an idealized public transportation system},
  author={Osuna, EE and Newell, Gordon F},
  journal={Transportation Science},
  volume={6},
  number={1},
  pages={52--72},
  year={1972},
  publisher={INFORMS}
}

@article{islam2015model,
  title={A model to evaluate the impact of headway variation and vehicle size on the reliability of public transit},
  author={Islam, Md Kamrul and Vandebona, Upali and Dixit, Vinayak V and Sharma, Ashish},
  journal={IEEE Transactions on Intelligent Transportation Systems},
  volume={16},
  number={4},
  pages={1840--1850},
  year={2015},
  publisher={IEEE}
}

@article{loynes1962stability,
  title={The stability of a queue with non-independent inter-arrival and service times},
  author={Loynes, R. M.},
  journal={Mathematical Proceedings of the Cambridge Philosophical Society},
  volume={58},
  number={3},
  pages={497--520},
  year={1962},
  doi={10.1017/S0305004100036781}
}

@article{lindley1952theory,
  title={The theory of queues with a single server},
  author={Lindley, D. V.},
  journal={Mathematical Proceedings of the Cambridge Philosophical Society},
  volume={48},
  number={2},
  pages={277--289},
  year={1952},
  doi={10.1017/S0305004100027638}
}
\endgroup

\ECRUNAUTHOR{}
\ECAUpunct{}
\ECSwitch
\providecommand{\theHsection}{}
\providecommand{\theHsubsection}{}
\providecommand{\theHsubsubsection}{}
\providecommand{\theHequation}{}
\providecommand{\theHfigure}{}
\providecommand{\theHtable}{}
\renewcommand{\theHsection}{EC.\arabic{section}}
\renewcommand{\theHsubsection}{EC.\arabic{section}.\arabic{subsection}}
\renewcommand{\theHsubsubsection}{EC.\arabic{section}.\arabic{subsection}.\arabic{subsubsection}}
\renewcommand{\theHequation}{EC.\arabic{equation}}
\renewcommand{\theHfigure}{EC.\arabic{figure}}
\renewcommand{\theHtable}{EC.\arabic{table}}
\ECHead{Electronic Companion}
\setcounter{algorithm}{0}
\renewcommand{\thealgorithm}{EC.\arabic{algorithm}}
\providecommand{\theHalgorithm}{}
\renewcommand{\theHalgorithm}{EC.\arabic{algorithm}}

\section{Modeling assumptions}\label{append_model_assumptions}

Table \ref{tab_ec_primitives} summarizes route scope and stochastic primitives, and Table \ref{tab_ec_dependence} summarizes dependence and service sequence approximations.

\begin{table}[htb]
\centering
\caption{Route scope and stochastic primitives}
\label{tab_ec_primitives}
\fontsize{8.6}{10.3}\selectfont
\begin{tabular}{@{}>{\raggedright\arraybackslash}p{0.27\textwidth}p{0.63\textwidth}@{}}
\toprule
Assumption & Operational interpretation \\ \midrule
Single route with ordered stations & The model targets line level resilience for one bus route or one rail direction, where vehicles serve stops in a fixed order. This is the natural unit for evaluating station queues, vehicle capacity propagation, and headway disruptions before adding transfer interactions. \\
Poisson passenger arrivals & Passenger arrivals are treated as random at the aggregate stop level during the analysis period. This is appropriate when individual arrivals are not perfectly synchronized with vehicle departures and when the arrival rate $\lambda^{(n)}$ can be estimated from observed demand. \\
Binomial alighting & Each onboard passenger has a station specific probability of alighting. This provides a parsimonious way to propagate vehicle load and downstream available capacity when only aggregate alighting shares or proportions by origin and destination are available. \\
Negligible dwell time & The analysis focuses on headway changes caused by short service suspensions rather than passenger service time at stops. When dwell time is small relative to the operating headway, passenger arrivals during dwell are second order, and systematic dwell components can be absorbed into the effective headway if needed. \\
No balking or reneging & The paper studies short suspensions, during which waiting passengers are more likely to remain in the system than abandon their trip. The assumption therefore keeps the resilience indicators focused on capacity, headway, and passengers left behind, without modeling route choice behavior during longer disruptions. \\
Full stop failure state & A full stop is a conservative first order representation of incidents such as signal failures, vehicle breakdowns, or blocked movements. Partial slowdowns can be approximated by an effective suspension duration or by extending the two state speed process. \\ \bottomrule
\end{tabular}
\end{table}

\begin{table}[htb]
\centering
\caption{Dependence and service sequence approximations}
\label{tab_ec_dependence}
\fontsize{8.6}{10.3}\selectfont
\begin{tabular}{@{}>{\raggedright\arraybackslash}p{0.27\textwidth}p{0.63\textwidth}@{}}
\toprule
Assumption & Operational interpretation \\ \midrule
No overtaking & Vehicles serve each station in dispatch order. This rule directly represents rail lines where passing is infeasible and bus services where passing is prohibited or uncommon. The recursive first in, first out (FIFO) model holds a vehicle when its free departure time precedes the realized departure time of the preceding vehicle. \\
Independent exogenous renewal headways & At each station, headways are independent across vehicle runs and independent of prior passenger queue and vehicle load states. Their common marginal law retains the disruption effect over one headway implied by two independent incident duration draws. Serial dependence and delay propagation across multiple vehicles are omitted. The approximation is most appropriate when holding or dispatch control resets the next service interval from the current departure rather than recovering the absolute timetable. It is assessed against a recursive FIFO simulation. \\
Marginal route decoupling & The assumption is exact at the first station. At downstream stations, separate local demand, independent renewal headways, and binomial alighting weaken dependence between the remaining onboard load and local queue when suspensions are short and stations are stable. The simulation assesses the approximation. \\
Exogenous station service sequence & Once the marginal available capacity distribution has been propagated from upstream, successive service opportunities at a station are independent draws from that distribution and are independent of the local queue and arrival processes. The complete capacity distribution carries the effect of upstream crowding into the station model, whereas the analytical queue omits temporal dependence in capacity. \\
Nonnegative headways without later recovery & Departure epochs are accumulated from nonnegative headways, so the analytical service sequence is nondecreasing by construction. A nonpositive raw headway becomes a zero headway, representing bunching, and the lost separation is not recovered in a later headway. The resulting analytical service sequence follows the renewal representation and differs from the realized FIFO trajectory. \\ \bottomrule
\end{tabular}
\end{table}

\clearpage

\section{Certified interpolation root search}\label{append_root_search}

Algorithm \ref{alg_inter_search} gives the complete procedure for the three search steps described in Section \ref{sec_root_solving} and reports the parameter settings used for every station and numerical scenario.

\begin{algorithm}[htb]
\caption{Certified interpolation search} \label{alg_inter_search}
\fontsize{9.5}{11.4}\selectfont
\begin{algorithmic}[1]
\Require Characteristic function $\mathcal{F}^{(n)}$, passenger arrival probability generating function $Y(z)$, capacity distribution $s^{(n)}$, capacity $C$, and utilization $\rho^{(n)}=\overline{Y}^{(n)}/\overline{S}^{(n)}$
\State Set the residual tolerance $\tau_{\mathrm{res}}\gets10^{-9}$ and clustering tolerance $\tau_{\mathrm{dup}}\gets10^{-4}$
\State Set the unit disk tolerance $\tau_{\mathrm{disk}}\gets10^{-8}$
\State Set $K_{\max}\gets\max\{12,C\}$ and $L_{\max}\gets\max\{4,\lceil C/2\rceil\}$
\State Define the scaled residual
\Statex \hspace{\algorithmicindent}$\displaystyle \operatorname{Res}^{(n)}(z):=\frac{|\mathcal{F}^{(n)}(z)|}{|z|^C+\left|Y(z)\sum_{u=0}^{C}s_u^{(n)}z^{C-u}\right|}$
\State For a candidate set $\mathcal{C}$, define $E_{\mathrm{dup}}(\mathcal{C}):=\{\{z,z'\}\subseteq\mathcal{C}:z\neq z',\ |z-z'|\leq\tau_{\mathrm{dup}}\}$
\Statex \hspace{\algorithmicindent}The connected components of $(\mathcal{C},E_{\mathrm{dup}}(\mathcal{C}))$ are the root clusters
\State Set $r_0\gets\min\{0.98,\max\{0.05,1-0.5\rho^{(n)}\}\}$
\State Initialize the candidate set $\mathcal{C}^{(0)}\gets\{1\}$
\For{$k=0,\ldots,C-1$}
    \State Set $z^{\mathrm{Ini}}\gets r_0\exp[2\pi(k+1/2)\sqrt{-1}/C]$ following \citet{powell1985analysis}
    \State Apply the local solver to $\mathcal{F}^{(n)}(z)=0$ from $z^{\mathrm{Ini}}$
    \If{the solver returns $z$ with $\operatorname{Res}^{(n)}(z)\leq\tau_{\mathrm{res}}$ and $|z|\leq1+\tau_{\mathrm{disk}}$}
        \State Add $z$ and its conjugate to $\mathcal{C}^{(0)}$
    \EndIf
\EndFor
\State Form the root clusters from $\mathcal{C}^{(0)}$ using the graph rule above and apply polar extrapolation
\State Certify the cluster multiplicities using Eq. \ref{eq_root_multiplicity}; obtain $\mathcal{Z}^{(0)}$, $M_0$, and the total certified multiplicity $N_0$
\State Set $L\gets2$ and $k\gets0$
\While{($N_k<C$ or at least one cluster is uncertified) and $k<K_{\max}$ and $L\leq L_{\max}$}
    \State Order $\mathcal{Z}^{(k)}=\{\zeta_0^{(k)},\ldots,\zeta_{M_k-1}^{(k)}\}$ by complex argument
    \State Initialize $\mathcal{Z}^{\mathrm{Ini}}\gets\emptyset$
    \For{$i=0,\ldots,M_k-1$}
        \State Set $j\gets(i+1)\bmod M_k$
        \State Set $r_i\gets|\zeta_i^{(k)}|$ and $r_j\gets|\zeta_j^{(k)}|$; unwrap $\phi_i\gets\arg(\zeta_i^{(k)})$ and $\phi_j\gets\arg(\zeta_j^{(k)})$ so that $\phi_j>\phi_i$
        \For{$d=1,\ldots,L-1$}
            \State Set $r_{i,d}\gets(1-d/L)r_i+(d/L)r_j$ and $\phi_{i,d}\gets(1-d/L)\phi_i+(d/L)\phi_j$
            \For{$a\in\{0.85,1,1.15\}$}
                \State Add $\min\{1,ar_{i,d}\}\exp(\phi_{i,d}\sqrt{-1})$ to $\mathcal{Z}^{\mathrm{Ini}}$
            \EndFor
        \EndFor
    \EndFor
    \State Apply the local solver from every initial value in $\mathcal{Z}^{\mathrm{Ini}}$ and refine every successful output
    \State Collect each output and its conjugate with $\mathcal{Z}^{(k)}$ in $\mathcal{C}^{(k+1)}$
    \State $\widetilde{\mathcal{C}}^{(k+1)}\gets\{z\in\mathcal{C}^{(k+1)}:\operatorname{Res}^{(n)}(z)\leq\tau_{\mathrm{res}},\ |z|\leq1+\tau_{\mathrm{disk}}\}$
    \State Let $\mathcal{C}_0^{(k+1)},\ldots,\mathcal{C}_{M_{k+1}-1}^{(k+1)}$ be the root clusters obtained from $\widetilde{\mathcal{C}}^{(k+1)}$ using the graph rule above
    \State Set $\zeta_j^{(k+1)}\gets|\mathcal{C}_j^{(k+1)}|^{-1}\sum_{z\in\mathcal{C}_j^{(k+1)}}z$ for each cluster; certify its multiplicity and obtain $\mathcal{Z}^{(k+1)}$, $M_{k+1}$, and $N_{k+1}$
    \State Remove an uncertified cluster only if the scaled residual of its center exceeds $\max\{10^{-13},\tau_{\mathrm{res}}/100\}$
    \If{$M_{k+1}=M_k$}
        \State Set $L\gets L+1$
    \EndIf
    \State Set $k\gets k+1$
\EndWhile
\If{$N_k=C$ and every cluster is certified}
    \State Return every cluster center repeated according to its certified multiplicity
\Else
    \State Return a diagnostic failure with $N_k$, the largest residual, $k$, and $L$
\EndIf
\end{algorithmic}
\end{algorithm}

\clearpage

\section{Proofs}

\subsection{Proof of Proposition \ref{prop_space_dist}}

Suppose that vehicle $l$ carries $i$ passengers when it reaches station $n$. During alighting, each passenger independently remains onboard with probability $1-\alpha^{(n)}$. The remaining load therefore follows a binomial distribution with parameters $i$ and $1-\alpha^{(n)}$. For every $j=0,\ldots,i$, the corresponding conditional probability is $a_{ij}^{(n)}=\binom{i}{j}(1-\alpha^{(n)})^j(\alpha^{(n)})^{i-j}$, as stated in Eq. \ref{eq_A_mat}.

The law of total probability gives $g_j^{(n,l)}=\sum_{i=0}^{C}v_i^{(n-1,l)}a_{ij}^{(n)}$. In vector form, this relation is
\begin{align}
g^{(n,l)} = v^{(n-1,l)}A^{(n)},
\label{eq_g1}
\end{align}
where $v^{(n-1,l)} = [v^{(n-1,l)}_0,\ldots,v^{(n-1,l)}_C]$ and $g^{(n,l)} = [g^{(n,l)}_0,\ldots,g^{(n,l)}_C]$. Equation \ref{eq_Sn_Gn} shows that $S^{(n,l)}=C-G^{(n,l)}$. Hence, the distribution of the available spaces after alighting satisfies
\begin{align}
s_k^{(n,l)} = g_{C-k}^{(n,l)} \quad\quad \forall k = 0,1,\ldots,C.
\label{eq_s1}
\end{align}

Equations \ref{eq_g1} and \ref{eq_s1} hold for every vehicle run. When the steady state limits defined in Section \ref{sec_space_dist} exist, taking $l\to\infty$ converts these run specific relations into their steady state counterparts and thereby establishes Proposition \ref{prop_space_dist}.

\subsection{Proof of Proposition \ref{prop_load_dist}}

After alighting and boarding at station $n$, the steady state vehicle load satisfies $V^{(n)}=\min\{C,G^{(n)}+Q^{(n)}\}$. Define $b_{ij}^{(n)}:=\mathbb{P}(V^{(n)}=j\mid G^{(n)}=i)$. If $0\leq i\leq j<C$, then $b_{ij}^{(n)}=\mathbb{P}(Q^{(n)}=j-i\mid G^{(n)}=i)$. By Assumption \ref{assumption_route_decoupling}, this conditional probability equals the marginal probability $q_{j-i}^{(n)}$, which gives the first case in Eq. \ref{eq_B_mat}.

If $j=C$ and $0\leq i<C$, the vehicle reaches capacity when $Q^{(n)}\geq C-i$. Assumption \ref{assumption_route_decoupling} therefore gives $b_{iC}^{(n)}=\mathbb{P}(Q^{(n)}\geq C-i)=1-\sum_{k=0}^{C-i-1}q_k^{(n)}$. If $i=j=C$, the vehicle is already full after alighting, so $b_{CC}^{(n)}=1$. All other transitions are infeasible and have probability zero.

The law of total probability now gives $v_j^{(n)}=\sum_{i=0}^{C}g_i^{(n)}b_{ij}^{(n)}$ for every $j=0,\ldots,C$. Writing these equations in matrix form yields $v^{(n)}=g^{(n)}B^{(n)}$, which completes the proof.

\subsection{Proof of Proposition \ref{prop_pgf_of_q}}\label{append_PGF_of_Q}

The proof adapts the transform argument of \citet{powell1981stochastic} to the arbitrary available capacity distribution $s_0^{(n)},\ldots,s_C^{(n)}$, consistent with the variable capacity transform of \citet{powell1985analysis}.

We begin with the one step transition of the queue. From the relationship in Eq. \ref{eq_relation_QRY},
\begin{align}
q^{(n,l+1)}_k = \sum_{i=0}^{k} r^{(n,l)}_i y_{k-i}^{(n,l)},
\label{eq_queueing_prob1}
\end{align}
where $r_k^{(n,l)} := \mathbb{P}(R^{(n,l)}=k)$ and $y_k^{(n,l)} := \mathbb{P}(Y^{(n,l)}=k)$ for every nonnegative integer $k$. For $u=0,\ldots,C$, also write $\left.r^{(n,l)}_k\right\vert_{S^{(n,l)}=u} := \mathbb{P}(R^{(n,l)}=k\mid S^{(n,l)}=u)$. Assumption \ref{assumption_station_service} makes $S^{(n,l)}$ independent of $Q^{(n,l)}$. Conditional on $S^{(n,l)}=u$, all passengers can board when $u\geq Q^{(n,l)}$, so no passenger is left behind. Therefore,
\begin{align}
\left.r^{(n,l)}_0\right\vert_{S^{(n,l)}=u}  = \mathbb{P}(u \geq Q^{(n,l)}\mid S^{(n,l)}=u) = \mathbb{P}(u \geq Q^{(n,l)}) = \sum_{k = 0}^u q^{(n,l)}_k.
\label{eq_relation_rq1}
\end{align}
If $u<Q^{(n,l)}$, only $u$ passengers can board and $Q^{(n,l)}-u$ passengers remain. Hence,

\begin{align}
\left.r^{(n,l)}_i\right\vert_{S^{(n,l)}=u}  = \mathbb{P}(Q^{(n,l)} - u = i\mid S^{(n,l)}=u) = \mathbb{P}(Q^{(n,l)} - u = i) = q^{(n,l)}_{u+i}, \qquad i=1,2,\ldots.
\label{eq_relation_rq2}
\end{align}

Using Eqs. \ref{eq_relation_rq1} and \ref{eq_relation_rq2}, Eq. \ref{eq_queueing_prob1} becomes
\begin{align}
q^{(n,l+1)}_k = \sum_{u=0}^C s^{(n,l)}_u \left(\sum_{i=0}^{u} q^{(n,l)}_i y_{k}^{(n,l)} + \sum_{i=u+1}^{u+k} q^{(n,l)}_i y_{k-i+u}^{(n,l)}\right).
\label{eq_queueing_prob2}
\end{align}

Assuming that the steady state probabilities exist, we have $\lim_{l\to\infty}q_k^{(n,l)}=q_k^{(n)}$, $\lim_{l\to\infty}s_k^{(n,l)}=s_k^{(n)}$, and $\lim_{l\to\infty}y_k^{(n,l)}=y_k^{(n)}$. Taking limits on both sides of Eq. \ref{eq_queueing_prob2} gives
\begin{align}
q^{(n)}_k = \sum_{u=0}^C s^{(n)}_u\sum_{i=0}^{u} q^{(n)}_i y_{k}^{(n)} + \sum_{u=0}^C s^{(n)}_u\sum_{i=u+1}^{u+k} q^{(n)}_i y_{k-i+u}^{(n)}.
\end{align}

Let $Q(z)$, $R(z)$, and $Y(z)$ denote the probability generating functions (PGFs) of the steady state random variables $Q^{(n)}$, $R^{(n)}$, and $Y^{(n)}$, respectively. Suppressing the station superscript, we write
\begin{align}
Q(z) &= \sum_{k=0}^{\infty} q^{(n)}_k z^k,\\
R(z) &= \sum_{k=0}^{\infty} r^{(n)}_k z^k = \sum_{k=0}^{\infty}  z^k \sum_{u=0}^C  s^{(n)}_u\cdot \left.r^{(n)}_k\right\vert_{S^{(n)}=u}, \label{eq_PGF_R}\\
Y(z) &= \sum_{k=0}^{\infty} y^{(n)}_k z^k.
\end{align}
Substituting Eqs. \ref{eq_relation_rq1} and \ref{eq_relation_rq2} into Eq. \ref{eq_PGF_R} gives
\begin{align}
R(z) &= \sum_{u=0}^C  s^{(n)}_u \sum_{i=0}^{u}q_i^{(n)} + \sum_{k=1}^{\infty}  z^k \sum_{u=0}^C s^{(n)}_u q_{k+u}^{(n)}\\
& = \sum_{u=0}^C s^{(n)}_u\left[\sum_{i=0}^{u}q_i^{(n)} + \frac{1}{z^u}Q(z) - \frac{1}{z^u}\sum_{i=0}^u q_i^{(n)} z^i \right].
\label{eq_R_intermsof_Q}
\end{align}

Taking the steady state limit in Eq. \ref{eq_relation_QRY} gives $Q^{(n)}=R^{(n)}+Y^{(n)}$. For a fixed vehicle run, the residual queue $R^{(n)}$ is determined by the queue and service outcomes through the current departure. Under Assumption \ref{assumption_renewal}, $Y^{(n)}$ is generated by the exogenous Poisson process during the next independently drawn headway and is therefore independent of this history. Consequently, $R^{(n)}$ and $Y^{(n)}$ are independent in the analytical model, and
\begin{align}
Q(z) = R(z)Y(z).
\label{eq_Q_eq_RY}
\end{align}
Combining Eqs. \ref{eq_R_intermsof_Q} and \ref{eq_Q_eq_RY} yields
\begin{align}
Q(z) &= \frac{Y(z)\sum_{u=0}^{C}s^{(n)}_u\left[\sum_{i=0}^u q_i^{(n)} (1-\frac{z^i}{z^u}) \right]}{1 - \sum_{u=0}^{C}s_u^{(n)} \frac{Y(z)}{z^u}} \nonumber\\
&= \frac{\sum_{u=0}^{C}s^{(n)}_u\left[\sum_{i=0}^u q_i^{(n)} (z^C-z^{C-u+i}) \right]}{\frac{z^C}{Y(z)} - \sum_{u=0}^{C}s_u^{(n)} z^{C-u}}.
\end{align}

\subsection{Proof of Proposition \ref{prop_solve_q}} \label{append_queue_dist}

Section \ref{sec_pgf_Q} identifies the $C$ characteristic roots of the denominator but does not determine the queue probabilities. We now use these roots to factor the numerator of $Q(z)$ and then match its polynomial coefficients to obtain $q_0^{(n)},\ldots,q_{C-1}^{(n)}$.

The numerator of $Q(z)$ is a polynomial of degree at most $C$. In Eq. \ref{eq_Qz1}, the term with $i=u$ vanishes, so the coefficient of $z^C$ in this numerator is $\kappa^{(n)}:=\sum_{u=1}^{C}s_u^{(n)}\sum_{i=0}^{u-1}q_i^{(n)}$. Analyticity requires the numerator to cancel all $C$ characteristic roots according to multiplicity. Therefore, the root $z_i^*$ is repeated in the following product whenever its multiplicity exceeds one, and $Q(z)$ can be written as
\begin{align}
Q(z) &= \frac{\kappa^{(n)}(z-1)\prod_{i=1}^{C-1}(z-z_i^*)}{\frac{z^C}{Y(z)} - \sum_{u=0}^{C}s_u^{(n)} z^{C-u}}.
\label{eq_Qz2}
\end{align}

As $z\to1$, both $\textsc{Num}(z)$ and $\textsc{Den}(z)$ approach zero, while $Q(z)\to1$. Applying l'H\^{o}pital's rule to this limit gives the following normalization identity.
\begin{align}
&\lim_{z\to1} Q(z) = 1 = \lim_{z\to1} \frac{\textsc{Num}'(z)}{\textsc{Den}'(z)} = \frac{\kappa^{(n)}\prod_{i=1}^{C-1}(1-z_i^*)}{\sum_{u=0}^{C}s_u^{(n)}u - Y'(1) } \\
\Rightarrow & \kappa^{(n)} = \frac{\sum_{u=0}^{C}s_u^{(n)}u - Y'(1) }{\prod_{i=1}^{C-1}(1-z_i^*)}. \label{eq_sq}
\end{align}

Recall that $\overline{Y}^{(n)} := Y'(1)=\mathbb{E}[Y^{(n)}]$ is the mean number of passenger arrivals within a headway and $\overline{S}^{(n)} := \sum_{u=0}^{C}s_u^{(n)}u=\mathbb{E}[S^{(n)}]$ is the mean number of available spaces in an arriving vehicle. Substituting Eq. \ref{eq_sq} into Eq. \ref{eq_Qz2} gives
\begin{align}
Q(z) &= \frac{(\overline{S}^{(n)} - \overline{Y}^{(n)})  (z-1)\prod_{i=1}^{C-1}\frac{z-z_i^*}{1-z_i^*}}{\frac{z^C}{Y(z)} - \sum_{u=0}^{C}s_u^{(n)} z^{C-u}}.
\label{eq_Qz3}
\end{align}

Equating the numerators in Eqs. \ref{eq_Qz3} and \ref{eq_Qz1} gives
\begin{align}
(\overline{S}^{(n)} - \overline{Y}^{(n)})  (z-1)\prod_{i=1}^{C-1}\frac{z-z_i^*}{1-z_i^*} = \sum_{u=0}^{C}s^{(n)}_u\left[\sum_{i=0}^u q_i^{(n)} (z^C-z^{C-u+i}) \right].
\label{eq_Qz_num}
\end{align}
Both sides of Eq. \ref{eq_Qz_num} are polynomials in $z$, so their coefficients must be equal. Matching the constant terms on both sides gives the first coefficient relation.
\begin{align}
(-1)(\overline{S}^{(n)} - \overline{Y}^{(n)})\prod_{i=1}^{C-1}\frac{z_i^*}{z_i^*-1} = -q_0^{(n)}s_C^{(n)}.
\end{align}
Solving this coefficient relation for the empty queue probability gives
\begin{align}
q_0^{(n)} = \frac{1}{s_C^{(n)}} (\overline{S}^{(n)} - \overline{Y}^{(n)})\prod_{i=1}^{C-1}\frac{z_i^*}{z_i^*-1}.
\end{align}

The coefficients of the characteristic equation are real, so every nonreal root and its multiplicity are paired with the corresponding conjugate root. The product in Eq. \ref{eq_q0_final} is therefore real. Under stability and irreducibility, the stationary empty queue probability is positive, which establishes the strict positivity of the product.

As a consistency check, consider fixed capacity with $s_C^{(n)}=1$ and $\overline{S}^{(n)}=C$. Under this specialization, Eq. \ref{eq_Qz_num} reduces to the following fixed capacity expression.
\begin{align}
\left.q_0^{(n)}\right\vert_{s_C^{(n)} = 1} = (C - \overline{Y}^{(n)})\prod_{i=1}^{C-1}\frac{z_i^*}{z_i^*-1}.
\end{align}
This expression agrees with \citet{chaudhry1987computational}. We next derive $q_{1:C-1}^{(n)}$. The numerator of Eq. \ref{eq_Qz3} can be rewritten as
\begin{align}
(\overline{S}^{(n)} - \overline{Y}^{(n)})  (z-1)\prod_{i=1}^{C-1}\frac{z-z_i^*}{1-z_i^*} &= \frac{1}{s_C^{(n)}}  (\overline{S}^{(n)} - \overline{Y}^{(n)}) \prod_{i=1}^{C-1}\frac{z_i^*}{z_i^*-1} \prod_{i=1}^{C-1}\frac{z_i^* - z}{z_i^*}  (z-1)s_C^{(n)} \nonumber \\
&= s_C^{(n)} q_0^{(n)} (z-1) \prod_{i=1}^{C-1} \left(1 - \frac{z}{z_i^*}\right)\nonumber \\
&=  - s_C^{(n)} q_0^{(n)} \prod_{i=0}^{C-1} \left(1 - \frac{z}{z_i^*}\right).
\label{eq_Qz_num_2}
\end{align}
Define $\prod_{i=0}^{C-1}(1-z/z_i^*) := \sum_{j=0}^{C}\eta_j^{(n)}z^j$, where $\eta_j^{(n)}$ is the coefficient of $z^j$ for station $n$. On the right hand side of Eq. \ref{eq_Qz_num}, the coefficient of $z^{C-k}$ is $-\sum_{u=k}^Cs_u^{(n)}q_{u-k}^{(n)}$. In Eq. \ref{eq_Qz_num_2}, the corresponding coefficient is $-s_C^{(n)}q_0^{(n)}\eta_{C-k}^{(n)}$. Matching these coefficients gives
\begin{align}
s_C^{(n)} q_0^{(n)} \eta_{C-k}^{(n)} = \sum_{u=k}^C s_u^{(n)} q_{u-k}^{(n)}, \qquad k=1,2,\ldots,C-1.
\label{eq_q_sq}
\end{align}
For fixed capacity, $s_C^{(n)}=1$ and $s_k^{(n)}=0$ for $0\leq k<C$. Equation \ref{eq_q_sq} then reduces to

\begin{align}
q_{C-k}^{(n)} = q_0^{(n)} \eta_{C-k}^{(n)}, \qquad k=1,2,\ldots,C-1, \quad \text{if }s_C^{(n)}=1,
\end{align}
which agrees with the result of \citet{chaudhry1987computational}.

Adding the identity $s_C^{(n)}q_0^{(n)}\eta_0^{(n)}=s_C^{(n)}q_0^{(n)}$, where $\eta_0^{(n)}=1$, converts Eq. \ref{eq_q_sq} into the following matrix representation used to solve the queue probabilities.
\begin{align}
\widetilde{\eta}^{(n)} = q_{0:C-1}^{(n)}\Lambda^{(n)},
\end{align}
where $\widetilde{\eta}^{(n)} := [s_C^{(n)}q_0^{(n)}\eta_0^{(n)},s_C^{(n)}q_0^{(n)}\eta_1^{(n)},\ldots,s_C^{(n)}q_0^{(n)}\eta_{C-1}^{(n)}]\in\mathbb{R}^{C}$ and

\begin{equation}
 \Lambda^{(n)} := \left[
    \begin{array}{ccccc}
     s_C^{(n)} & s_{C-1}^{(n)} & \cdots & s_2^{(n)} & s_1^{(n)} \\
     0 & s_C^{(n)} & \cdots & s_3^{(n)} & s_2^{(n)} \\
     \vdots & \ddots & \ddots & \vdots & \vdots \\
     0 & \cdots & 0 & s_C^{(n)} & s_{C-1}^{(n)} \\
     0 & \cdots & \cdots & 0 & s_C^{(n)}
    \end{array}\right] \in \mathbb{R}^{C\times C}.
\end{equation}
Because $s_{C}^{(n)} > 0$, the triangular matrix $\Lambda^{(n)}$ is invertible. Therefore,
\begin{align}
q_{0:C-1}^{(n)} = \widetilde{\eta}^{(n)} (\Lambda^{(n)})^{-1}.
\end{align}

\subsection{Proof of Proposition \ref{prop_q_length}} \label{append_queue}

The derivation adapts the transform argument of \citet{powell1981stochastic} and is consistent with the general bulk service result of \citet{powell1985analysis}. For $X\in\{S,Y\}$, the symbols $\overline{X}^{(n)}$, $\overline{\overline{X}}^{(n)}$, and $\overline{\overline{\overline{X}}}^{(n)}$ denote its mean, second central moment, and third central moment, respectively. Define
\begin{align}
\mathcal{A}(z) &:= \frac{(\overline{S}^{(n)}-\overline{Y}^{(n)})(z-1)}{\frac{z^C}{Y(z)}-\sum_{u=0}^{C}s_u^{(n)}z^{C-u}},
& \mathcal{B}_i(z) &:= \frac{z-z_i^*}{1-z_i^*}.
\label{eq_AB_definition}
\end{align}
Equation \ref{eq_Qz3} gives $Q(z)=\mathcal{A}(z)\prod_{i=1}^{C-1}\mathcal{B}_i(z)$. At $z=1$, $Q(1)=\mathcal{A}(1)=\mathcal{B}_i(1)=1$ and $\mathcal{B}_i'(1)=1/(1-z_i^*)$. Differentiating the product then gives $Q'(1)=\mathcal{A}'(1)+\sum_{i=1}^{C-1}1/(1-z_i^*)$.

Let $\mathcal{A}_1(z):=(\overline{S}^{(n)}-\overline{Y}^{(n)})(z-1)$ and let $\mathcal{A}_2(z)$ denote the denominator of $\mathcal{A}(z)$ in Eq. \ref{eq_AB_definition}. We have $\mathcal{A}_1(1)=\mathcal{A}_2(1)=0$ and $\mathcal{A}_1'(1)=\mathcal{A}_2'(1)=\overline{S}^{(n)}-\overline{Y}^{(n)}>0$. A Taylor expansion of $\mathcal{A}_2(z)$ around $z=1$ gives the following two derivative identities.
\begin{align}
\mathcal{A}'(1) &= -\frac{\mathcal{A}_2''(1)}{2\mathcal{A}_2'(1)},
& \mathcal{A}''(1) &= \frac{3[\mathcal{A}_2''(1)]^2-2\mathcal{A}_2'(1)\mathcal{A}_2'''(1)}{6[\mathcal{A}_2'(1)]^2}.
\label{eq_A_derivatives}
\end{align}
Differentiating $\mathcal{A}_2(z)$ and expressing the resulting raw moments in terms of central moments gives
\begin{align}
\mathcal{A}'(1)=\frac{\overline{\overline{S}}^{(n)}+\overline{\overline{Y}}^{(n)}+(\overline{S}^{(n)}-\overline{Y}^{(n)})[1+2(\overline{S}^{(n)}-C)]-(\overline{S}^{(n)}-\overline{Y}^{(n)})^2}{2(\overline{S}^{(n)}-\overline{Y}^{(n)})}.
\end{align}
Substitution into $Q'(1)$ yields the mean in Eq. \ref{eq_EQ}.

For the variance, differentiating $\log Q(z)=\log \mathcal{A}(z)+\sum_{i=1}^{C-1}\log \mathcal{B}_i(z)$ twice and then applying $\mathcal{A}(1)=\mathcal{B}_i(1)=1$ and $\mathcal{B}_i''(1)=0$ yields the following identity.
\begin{align}
\operatorname{Var}[Q^{(n)}] = \mathcal{A}''(1) - \mathcal{A}'(1)^2  + \mathcal{A}'(1) + \sum_{i=1}^{C-1} \left(\mathcal{B}_i'(1) - [\mathcal{B}_i'(1)]^2 \right).
\label{eq_var_Q2}
\end{align}

Using Eq. \ref{eq_A_derivatives}, differentiating $\mathcal{A}_2(z)$ through third order, and then converting the resulting raw moments to central moments yields the following expression.
\begin{align}
&\mathcal{A}''(1)-\mathcal{A}'(1)^2+\mathcal{A}'(1) \nonumber\\
&=\frac{1}{12(\overline{S}^{(n)}-\overline{Y}^{(n)})^2}\bigg[
-4(\overline{\overline{\overline{S}}}^{(n)}-\overline{\overline{\overline{Y}}}^{(n)})(\overline{S}^{(n)}-\overline{Y}^{(n)})
+3(\overline{\overline{S}}^{(n)}+\overline{\overline{Y}}^{(n)})^2 \nonumber\\
&\qquad-[6(\overline{\overline{S}}^{(n)}-\overline{\overline{Y}}^{(n)})-1](\overline{S}^{(n)}-\overline{Y}^{(n)})^2
-(\overline{S}^{(n)}-\overline{Y}^{(n)})^4\bigg].
\label{eq_Q_var_first_part}
\end{align}
Finally, $\mathcal{B}_i'(1)-[\mathcal{B}_i'(1)]^2=-z_i^*/(1-z_i^*)^2$. Substituting this identity and Eq. \ref{eq_Q_var_first_part} into Eq. \ref{eq_var_Q2} yields the variance in Proposition \ref{prop_q_length}.

\subsection{Proof of Proposition \ref{prop_hdw_distribution}} 

Without incidents, vehicle $l$ requires $T^{(n)}$ time units to reach station $n$. Incident onsets occur only during normal travel and follow a Poisson process with rate $\gamma$, so the number of incidents during this travel follows a Poisson distribution with mean $\gamma T^{(n)}$. Independent exponential incident durations with rate $\theta$ then yield the compound Poisson representation in Eq. \ref{eq_I_sum}.

\subsection{Proof of Proposition \ref{prop_headway_distribution1}}

The raw headway moment generating function (MGF) is
\begin{align}
 M_{\widehat{H}^{(n,l)}}(t) &= \mathbb{E}[e^{t \widehat{H}^{(n,l)}}]
 = \mathbb{E}\left[e^{t(\overline{H} + \frac{2 \gamma T^{(N)}}{\theta\overline{F}} + I^{(n,l)} - I^{(n,l-1)})}\right] \notag\\
 & = e^{t(\overline{H} + \frac{2 \gamma T^{(N)}}{\theta\overline{F}})}
 \mathbb{E}[e^{tI^{(n,l)}}]  \mathbb{E}[e^{-tI^{(n,l-1)}}] \notag\\
& = e^{t(\overline{H} + \frac{2 \gamma T^{(N)}}{\theta\overline{F}})}  e^{\gamma T^{(n)} (\frac{\theta}{\theta - t} - 1)}  e^{\gamma T^{(n)} (\frac{\theta}{\theta + t} - 1)} \notag\\
& =  e^{t(\overline{H} + \frac{2 \gamma T^{(N)}}{\theta\overline{F}})}  e^{\gamma T^{(n)} (\frac{2t^2}{\theta^2 - t^2})}, \qquad |t|<\theta. \label{eq_H_mgf_ind}
\end{align}
The factorization in Eq. \ref{eq_H_mgf_ind} follows from the independence of the two copies $I^{(n,l)}$ and $I^{(n,l-1)}$ used to define one generic raw headway. Because this is the common marginal law used for every renewal draw, $M_{\widehat{H}^{(n)}}(t)=M_{\widehat{H}^{(n,l)}}(t)$. Assumption \ref{assumption_renewal} separately imposes independence across successive headway draws in the analytical renewal sequence.

\subsection{Proof of Proposition \ref{prop_headway_distribution}}

Let $\mu = \overline{H} + \frac{2 \gamma T^{(N)}}{\theta\overline{F}}$ and $\sigma = \frac{2 \sqrt{T^{(n)} \gamma}}{\theta}$ denote the mean and standard deviation of $\widehat{H}^{(n)}_{\text{Normal}}$, respectively. Because $H^{(n)}_{\text{Normal}}$ is the positive part of this variable, its MGF is
\begin{align}
  M_{H^{(n)}_{\text{Normal}}}(t) &= \mathbb{E}[e^{t H^{(n)}_{\text{Normal}}}] = \mathbb{P}[\widehat{H}^{(n)}_{\text{Normal}}\leq0]\cdot e^{0} + \int_{0}^{+\infty}  e^{t z} \cdot \phi_{\widehat{H}^{(n)}_{\text{Normal}}}{(z)} \cdot dz \notag\\
  &= \Phi(\frac{-\mu}{\sigma}) + \frac{1}{\sigma\sqrt{2\pi}} \int_{0}^{+\infty}  e^{t z + \frac{(z - \mu)^2}{-2\sigma^2}} dz \notag\\
  & =  \Phi(\frac{-\mu}{\sigma}) + e^{\mu t + \frac{\sigma^2t^2}{2}}\left[1 - \Phi(\frac{-\mu}{\sigma} - \sigma t)\right]. \label{eq_positive_normal_mgf}
%   &= \Phi\left(\frac{-(\overline{H}\theta + \frac{2\gamma T^{(N)}}{\overline{F}})}{2\sqrt{T^{(n)}\gamma}}\right) +  e^{t(\overline{H} + \frac{2 \gamma T^{(N)}}{\theta\overline{F}})}  e^{\gamma T^{(n)}(\frac{2t^2}{\theta^2})}\left[1 -  \Phi\left(\frac{-(\overline{H}\theta + \frac{2\gamma T^{(N)}}{\overline{F}})}{2\sqrt{T^{(n)}\gamma}} - \frac{2t\sqrt{T^{(n)} \gamma}}{\theta}\right)\right]
\end{align}
Completing the square in the exponent gives the final equality in Eq. \ref{eq_positive_normal_mgf}. Substituting the definitions of $\mu$ and $\sigma$ into this expression establishes Proposition \ref{prop_headway_distribution}.

\subsection{Proof of Proposition \ref{prop_headway_increase}}

\begin{figure}[htb]
\centering
\subfloat{\includegraphics[width=0.5\textwidth]{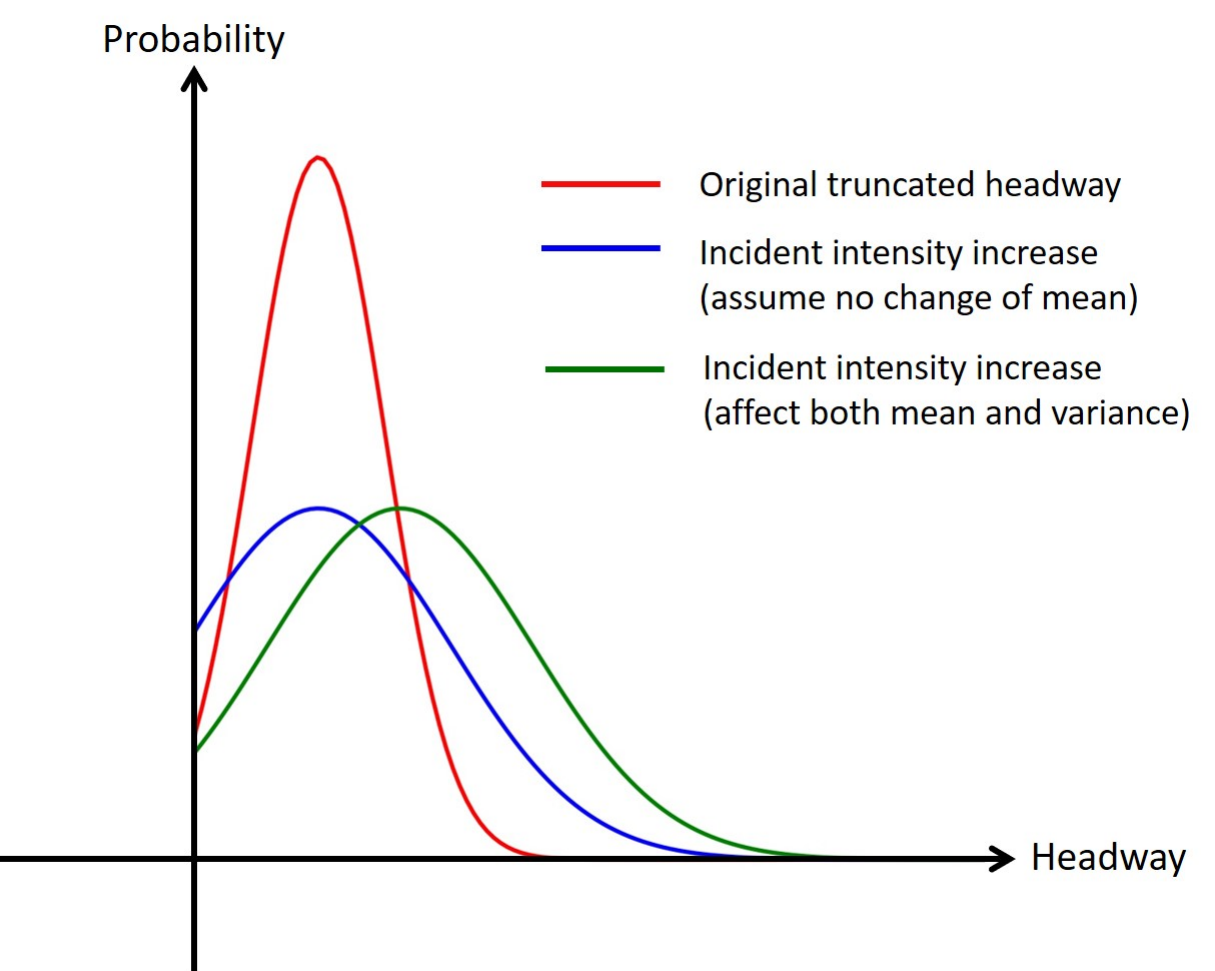}}
\caption{Illustration of the impact of incidents on expected headway. As the probability mass at zero does not contribute to the expectation calculation, it is not shown in the figure.}
\label{fig_prove_hdw_inc}
\end{figure}

For a normal variable with mean $\mu$ and standard deviation $\sigma>0$, Eq. \ref{eq_mean_normal_tr_hdw} gives the expectation of its positive part. Differentiating this expectation with respect to its two parameters gives
\begin{align}
\frac{\partial\mathbb{E}[H^{(n)}_{\text{Normal}}]}{\partial\mu}
&=\Phi\left(\frac{\mu}{\sigma}\right)>0,
&
\frac{\partial\mathbb{E}[H^{(n)}_{\text{Normal}}]}{\partial\sigma}
&=\phi\left(\frac{\mu}{\sigma}\right)>0.
\label{eq_positive_normal_derivatives}
\end{align}
The mean $\mu=\overline{H}+2\gamma T^{(N)}/(\theta\overline{F})$ and standard deviation $\sigma=2\sqrt{T^{(n)}\gamma}/\theta$ are nondecreasing in both $\gamma$ and $1/\theta$. The mean is strictly increasing in $\gamma$, and both parameters are strictly increasing in $1/\theta$ when $\gamma>0$. Equation \ref{eq_positive_normal_derivatives} establishes the stated monotonicity. The case $\gamma=0$ follows by continuity. Figure \ref{fig_prove_hdw_inc} illustrates these two effects.

\subsection{Proof of Proposition \ref{prop_arrival_pax}}

Conditional on $H^{(n)}_{\text{Normal}}$, the variable $Y^{(n)}$ is Poisson with parameter $\lambda^{(n)}H^{(n)}_{\text{Normal}}$. The law of iterated expectations therefore gives $Y(z)=\mathbb{E}\{\exp[\lambda^{(n)}H^{(n)}_{\text{Normal}}(z-1)]\}=M_{H^{(n)}_{\text{Normal}}}[\lambda^{(n)}(z-1)]$. For $\gamma>0$, substituting $t=\lambda^{(n)}(z-1)$ into the MGF in Proposition \ref{prop_headway_distribution} gives the stated PGF. When $\gamma=0$, $H^{(n)}_{\text{Normal}}=\overline{H}$ almost surely, so the conditional Poisson PGF directly gives $Y(z)=\exp[\lambda^{(n)}\overline{H}(z-1)]$.

\subsection{Proof of Proposition \ref{prop_stability}}

We establish stability directly from the transient queue dynamics. By the definitions of the queue and available capacity, $R^{(n,l)}=\max\{0,Q^{(n,l)}-S^{(n,l)}\}$. Combining this relation with Eq. \ref{eq_relation_QRY} gives
\begin{align}
R^{(n,l+1)}=\max\{0,R^{(n,l)}+Y^{(n,l)}-S^{(n,l+1)}\}.
\label{eq_lindley_residual}
\end{align}
Following \citet{lindley1952theory}, Eq. \ref{eq_lindley_residual} represents the residual passenger queue after successive vehicle service events as the standard Lindley recursion.

Assumption \ref{assumption_renewal}, the independent increments of the Poisson passenger arrival process, and Assumption \ref{assumption_station_service} imply that the net inputs $Y^{(n,l)}-S^{(n,l+1)}$ are identically and independently distributed across vehicle runs. They are exogenous to the current residual queue and have finite mean $\overline{Y}^{(n)}-\overline{S}^{(n)}$. The stability theorem of \citet{loynes1962stability} therefore applies directly to Eq. \ref{eq_lindley_residual}. If $\overline{Y}^{(n)}<\overline{S}^{(n)}$, the recursion has a proper stationary solution. Because $Q^{(n,l+1)}=R^{(n,l)}+Y^{(n,l)}$, the queue observed at vehicle arrival also has a proper stationary distribution.

If $\overline{Y}^{(n)}>\overline{S}^{(n)}$, the strong law of large numbers gives a positive long run drift for the partial sums of the net inputs, and the residual queue grows without bound. At the boundary $\overline{Y}^{(n)}=\overline{S}^{(n)}$, the nonconstant independent net inputs in Assumption \ref{assumption_station_service} define a critical Lindley recursion. The critical case has no proper stationary probability distribution under the same result of \citet{loynes1962stability}. Hence, the station is stable if and only if $\overline{Y}^{(n)}<\overline{S}^{(n)}$. Dividing both sides by $\overline{S}^{(n)}>0$ and substituting the expressions for the two means gives Eq. \ref{eq_stab}.

Upstream marginal propagation supplies $s^{(n)}$, and the proof requires only that the resulting station service sequence satisfy Assumption \ref{assumption_station_service}. The assumption makes this sequence exogenous to the local queue history. The station by station argument therefore does not presuppose the existence of $q_0^{(n)}$ or another stationary queue probability at the station under consideration.

\section{Simulation procedure for comparison}\label{append_sim_alg}
The simulation deliberately relaxes the analytical renewal, marginal route decoupling, and exogenous station service sequence assumptions. For each vehicle, we generate independent compound Poisson exponential incident durations on successive route segments and cumulatively sum them to obtain $I^{(n,l)}$ at each station. Realized departure times follow the recursive FIFO rule in Eq. \ref{eq_fifo_departure}. Vehicle $l$ cannot depart station $n$ before vehicle $l-1$, so an extreme delay can propagate through several successors. Passenger queues and vehicle loads are updated jointly rather than through independent marginal distributions. When a vehicle arrives, passengers board according to the first come, first served principle up to capacity $C$. Queue lengths at vehicle arrival and waiting times of boarded passengers are recorded during the simulation.

We evaluate the reference setting and six variations that cover low capacity, long headway, high demand, no incidents, frequent incidents, and long incident duration. Each setting uses 20 independent random seeds. A replication contains 10,000 vehicle runs, and the first 1,000 runs are excluded to reduce initialization effects. We calculate each performance measure separately within a replication and then use variation across replications to construct a 95\% confidence interval. The relative $\ell_1$ errors compare the mean across replications with the analytical result over stations 1 through 9. Station 10 is excluded from the waiting time evaluation because its zero arrival rate produces no waiting time observations.
\begin{algorithm}[htb]
\caption{Simulation procedure}
\fontsize{9.5}{11.4}\selectfont
\begin{algorithmic}[1]
\For{each of the seven assessment settings}
    \For{each of 20 independent random seeds}
        \State Initialize parameters and empty passenger queues and vehicle loads
        \For{$l=1,\ldots,10{,}000$}
            \State Set the scheduled dispatch time $DT^{(l)}\gets(l-1)\overline{H}^{\mathrm{Adj}}$
            \For{$n=1,\ldots,N$}
                \State Sample the incident duration on segment $n$ and cumulatively obtain $I^{(n,l)}$
                \State $\hat{t}^{(n,l)}_D \gets DT^{(l)}+T^{(n)}+I^{(n,l)}$
                \If{$l=1$}
                    \State $t^{(n,l)}_D\gets\hat{t}^{(n,l)}_D$ and initialize the station clock
                \Else
                    \State $t^{(n,l)}_D \gets \max\{\hat{t}^{(n,l)}_D,t^{(n,l-1)}_D\}$
                    \State Sample passenger arrival times during $(t^{(n,l-1)}_D,t^{(n,l)}_D]$ from the station Poisson process
                \EndIf
                \State Record queue length, including passengers left behind by the preceding vehicle
                \State Sample alighting using $\alpha^{(n)}$ and board passengers by the first come, first served rule up to capacity $C$
                \State Record the residual queue and waiting times of boarded passengers
            \EndFor
        \EndFor
        \State Exclude the first 1,000 vehicle runs and calculate the four station performance measures
    \EndFor
    \State Calculate means and 95\% confidence intervals across replications and relative $\ell_1$ errors over stations 1 through 9
\EndFor
\end{algorithmic}
\end{algorithm}

%%%%%%%%%%%%%%%%%
\end{document}